\documentclass[12pt]{amsart}
\usepackage{amssymb}
\usepackage[dvips]{graphics}
\ifx\mathrm\undefined\let\mathrm\rm\fi
\ifx\mathbf\undefined\let\mathbf\bf\fi
\ifx\mathfrak\undefined\let\mathfrak\frak\fi
\ifx\mathcal\undefined\let\mathcal\cal\fi
\ifx\mathbb\undefined\let\mathbb\Bbb\fi
\ifx\emph\undefined\let\emph\it\fi

\newcommand{\h}{{{\mathfrak h\,}}}
\newcommand{\Id}{{\operatorname{Id}}}
\newcommand{\Z}{{\mathbb Z}}

\newcommand{\C}{{\mathbb C}}

\newcommand{\Ref}[1]{{(\ref{#1})}}
\newcommand{\be}{\begin{displaymath}}
\newcommand{\ee}{\end{displaymath}}
\newcommand{\bea}{\begin{eqnarray*}}
\newcommand{\eea}{\end{eqnarray*}}
\newcommand{\tr}{{\mathrm{tr}}}

\newcommand{\dontprint}[1]{\relax}
\renewcommand{\th}{{{}^\mathrm{th}}}
\newenvironment{prf}{\noindent{\it Proof\/}:}{$\;\square$
\par\medskip}

\newtheorem%
{thm}{Theorem}[section]
\newtheorem%
{proposition}[thm]{Proposition}
\newtheorem%
{lemma}[thm]{Lemma}
\newtheorem%
{lemmadef}[thm]{Lemma-Definition}
\newtheorem%
{corollary}[thm]{Corollary}
\newtheorem%
{conjecture}[thm]{Conjecture}

\title[{}]{Separation of variables for quantum integrable systems
on elliptic curves}
\author[{}]
{ Giovanni Felder \and Anke Schorr}
\address{
Departement Mathematik, ETH-Zentrum, CH-8092
Z\"urich, Switzerland}
\begin{document}
\begin{abstract}
We extend Sklyanin's method of separation of variables to quantum
integrable models associated to elliptic curves. After reviewing
the differential case, the elliptic Gaudin model studied by 
Enriquez, Feigin and Rubtsov, we consider the difference case and
find a class of transfer matrices whose eigenvalue problem can
be solved by separation of variables. These transfer matrices
are associated to representations of the elliptic quantum group
$E_{\tau,\eta}(sl_2)$ by difference operators.
One model of statistical mechanics to which this method applies
is the IRF model with antiperiodic boundary conditions. The 
eigenvalues of the transfer matrix 
are given as solutions of a system of quadratic equations in
a space of higher order theta functions.
\end{abstract}

\maketitle
\date{}
\section{Introduction}

The method of separation of variables in integrable lattice models,
proposed by Sklyanin, is a method to find eigenvalues and
eigenvectors of transfer matrices. It is an alternative to
the Bethe ansatz and works in some situations where the Bethe
ansatz does not, and gives an insight in the completeness
of the Bethe ansatz. The method is closely related
to Baxter's method (Chapter 9 of \cite{B}), and
in fact the eigenvalue problem in the separated variables
(in the difference case) becomes the Baxter difference equation.
In the Gaudin model, one of the simplest quantum
integrable systems, this method relates the problem of finding
common eigenvectors of Gaudin Hamiltonians to the problem of
finding differential equations on the Riemann sphere
with regular singular points whose monodromy is trivial. As 
noticed by Feigin and Frenkel (see \cite{F}), this is a
special case of the Beilinson-Drinfeld
``geometric Langlands correspondance'' relating
certain local systems on a complex curve
to $\mathcal{D}$-modules on moduli spaces of principal bundles
on the curve.

Both for quantum integrable models and
for the connection to the Langlands program, it may be interesting
to extend the method of separation of variables to more general
models.

The class of quantum integrable systems (families of commuting
operators) one considers in this context
arise in different classes. There are the differential
models, such as the quantization of the Hitchin systems. They
are given by families of commuting differential operators and
are associated to complex curves. The Gaudin operators are the
operators associated to the Riemann sphere. More generally one
considers difference or $q$-deformed models, such as the six-vertex model,
which degenerate to the differential model when the parameter
$q$ tends to one. These models appear in three sorts: rational,
trigonometric and elliptic, depending on the type of coefficients
in the commuting operators.

From the point of view of representation theory, the differential
models are related to Kac--Moody Lie algebras, and the $q$-deformed
models to (infinite dimensional) quantum groups. 

We consider here models related to $sl_2$. In the differential 
rational (genus zero) case
the separation of variables was considered by Sklyanin
\cite{S} and Frenkel \cite{F}. A  version
of the separation of variables for the genus one differential case
was considered by Enriquez, Feigin and Rubtsov \cite{EFR}, who
also made an explicit connection to the Langlands correspondance. The
$q$-deformed rational and trigonometric case were studied by
Sklyanin \cite{S} and Tarasov--Varchenko \cite{TV}. The latter
authors introduce the notion of a difference equation with regular
singular points, thus giving a $q$-version of the relation described
above for Gaudin models. 

Here we consider the $q$-deformed elliptic case.  The class
of difference operator we give is both a $q$-deformation of the
Enriquez--Feigin--Rubtsov differential operators and an elliptic
version of the operators studied by Tarasov--Varchenko. Common
eigenfuncions may be constructed by the method of separation of
variables. Moreover the commuting difference operators we introduce
can be restricted to functions on a finite subset of points. This
restriction turns out to give the commuting transfer matrices
of interaction-round-a-face models with antiperiodic boundary
condition. This model provides an example of a model solvable
by separation of variables but not by Bethe ansatz. Other
elliptic models, related to the XYZ model,
have recently been studied by the method
of separation of variables by Sklyanin and Takebe \cite{ST}.

The algebraic structure at the origin of our constructions is
the elliptic quantum group $E_{\tau,\eta}(sl_2)$. Indeed, the
starting point is the construction of a representation of
this quantum group by difference operators which generalizes
the ``universal evaluation module'' of \cite{FV3}.

The paper is organized as follows:

In Section 2 we review the separation of variables in the 
differential elliptic case. Most of this part is essentially taken from 
\cite{EFR}, but we add some remarks on the Bethe ansatz and
its completeness.

In Section 3 we explain what is needed from the theory of
elliptic quantum groups and introduce a class of representations
of $E_{\tau\eta}(sl_2)$ by difference operators and relate
them to known representations. Twisted
commuting transfer matrices are then introduced and the method
of separation of variables is applied to construct (Bethe ansatz)
eigenvectors. 

In Section 4 we consider the restriction of the transfer matrix
associated to the tensor product of $n$ fundamental representations
to functions of a finite set of cardinality $2^n$, and
show that we obtain the transfer matrix of an IRF model with
antiperiodic boundary conditions. The eigenvalues are then
obtained as the solutions of a system of $n$ quadratic equations
in an $n$-dimensional space of theta functions of order $n$.

In an appendix we give an account on ``elliptic polynomials'', which
are (twisted) theta functions of order $n$.

\section{The differential case}\label{s-dc}
Let us start by introducing a family of commuting
differential operators associated to an elliptic curve
with $n$ marked points and $n$ highest weight representations of 
$sl_2(\C)$. 

Let the elliptic curve be $E=\C/\Gamma$ with $\Gamma=
\Z+\tau\Z$ and Im$\,\tau>0$. The marked points
are the projections of $n$ points $z_1,\dots,z_n\in\C$,
with $z_i\neq z_j \mod \Gamma$ for $i\neq j$. 
The representations are $M_{\Lambda_1},\dots,M_{\Lambda_n}$
where, for $\Lambda\in\C$,
$M_\Lambda$ denotes the Verma module (defined below in \ref{ss-32}
of highest weight $\Lambda$ 
of $sl_2(\C)$.

Thus our parameters are $\tau,z_1,\dots,z_n,\Lambda_1,\dots,\Lambda_n$.

Let $\theta(z)=-\sum_{j\in\Z}\exp({i\pi (j+1/2)^2\tau
+2\pi i(j+1/2)(z+1/2)})$ be the odd Jacobi theta function,
and set $\sigma_\lambda(z)=\frac{\theta(\lambda-z)\theta'(0)}
{\theta(z)\theta(\lambda)}$. It is the unique meromorphic
function of $z$ regular on $\C-\Gamma$, with a simple
pole with residue one at 0, and such that $\sigma_\lambda(
z+r+s\tau)=e^{2\pi is\lambda}\sigma_\lambda(z)$, $r,s\in\Z$.

\newcommand{\mat}[4]
{\left(\begin{array}{cc}#1&#2\\ #3&#4\end{array}\right)}
Let $e=({}^{0\;1}_{0\;0})$, $f=({}^{0\;0}_{1\;0})$, $h=
({}^{1\;\;0}_{0\,-\!1})$ be the standard
generators of $sl_2(\C)$. For $a\in sl_2(\C)$ let 
$a^{(i)}$ denote the action of $a$ on the $i$th factor of the
tensor product $M=M_{\Lambda_1}\otimes \cdots\otimes M_{\Lambda_n}$.

Introduce the following endomorphisms of $M$ depending on
$z,\lambda\in\C$:
\[
h(z)=\sum_{i=1}^n
\frac{\theta'(z-z_i)}{\theta(z-z_i)}h^{(i)},
\qquad
e_\lambda(z)=\sum_{i=1}^n\sigma_{-\lambda}(z-z_i)e^{(i)},
\qquad
f_\lambda(z)=\sum_{i=1}^n\sigma_{\lambda}(z-z_i)f^{(i)}.
\]
The family of commuting differential operators
is then obtained by the following generating function,
which is an elliptic version of the generating function
of Gaudin Hamiltonians \cite{S}.

\begin{thm} Let for $z\in\C$,
$S(z)$  be the differential operator acting on functions
of one complex variables $\lambda$ with values in the
zero weight space $M[0]=\{v\in M,\sum_i h^{(i)}v=0\}$ of $M$.
\[
S(z)=\left(\frac{\partial}{\partial\lambda}-\frac12
h(z)\right)^2+e_\lambda(z)f_\lambda(z)+f_\lambda(z)e_\lambda(z).
\]
Then $S(z)S(w)=S(w)S(z)$.
\end{thm}

One way of proving this theorem is to notice that it is
a special case of the flatness of the KZB connection
(Prop.\ 2 in \cite{FV1}). The relation to the KZB connection
is the following.
In the $sl_2(\C)$ case, the KZB connection involves the
differential 
operators (appearing on the right hand side of the 
KZB equations)
\[
H_j=-h^{(j)}\frac\partial{\partial\lambda}
+\sum_{k:k\neq j}
\frac12
\frac{\theta'(z_j-z_k)}
{\theta(z_j-z_k)}
h^{(j)}h^{(k)}
+\sigma_\lambda(z_j-z_k)e^{(j)}f^{(k)}
+\sigma_{-\lambda}(z_j-z_k)f^{(j)}e^{(k)},
\]
$j=1,\dots,n$ and
\[
H_0=
\frac{\partial^2}{\partial\lambda^2}
+\frac12
\sum_{j,k=1}^n
\frac12h^{(j)}h^{(k)}
\frac{\theta''(z_j-z_k)}
{\theta(z_j-z_k)}
-(e^{(j)}f^{(k)}
+f^{(j)}e^{(k)})
\frac{\partial\sigma_\lambda(z_j-z_k)}
{\partial\lambda}.
\]
The fact that the connection is flat means in particular that
these operators
form a commuting family when acting on $M[0]$-valued
functions. The terms with $j=k$ are understood as the limit
as the argument $z_j-z_k$ tends to zero. 
So $\theta''(0)/\theta(0)$ means 
$\theta'''(0)/\theta'(0)$ and 
$(\partial\sigma_\lambda/\partial\lambda)(0)=(\theta'/\theta)'(\lambda)$.
 Note that $\sum_{j=1}^n H_j$ vanishes
on $M[0]$-valued functions. We call the operators $H_j$, $j\geq 1$, the
elliptic Gaudin Hamiltonian and $H_0$ the (generalized)
Lam\'e Hamiltonian (it is the Lam\'e operator if $n=1$).

The relation between these commuting operators and $S(z)$ is:

\begin{proposition}
        Let\footnote{The relation of these functions with
the classical Weierstrass $\zeta$ and $\wp$ functions $\zeta(z)=
\frac1z+\sum_{(r,s)\in \Z^2-(0,0)}
\frac1{z+r+s\tau}
-\frac1{r+s\tau}
+\frac z{(r+s\tau)^2}$, $\wp(z)=-\zeta'(z)$ is
$\zeta(z)=\bar\zeta(z)+2\eta_1 z$, $\wp(z)=\bar\wp(z)-2\eta_1$, where
$2\eta_1=\theta'''(0)/3\theta'(0)$.}
        $\bar\zeta(z)=\theta'(z)/\theta(z)$,
        $\bar\wp(z)=-\bar\zeta'(z)$. Then
        \[
          S(z)=\sum_{k=1}^n\frac{c_k}2\bar\wp(z-z_k)
              +\sum_{k=1}^nH_k\bar\zeta(z-z_k)
              +H_0,
        \]
and $c_j$ is the Casimir value $c_j=\frac12\Lambda_j(\Lambda_j+2)$.
\end{proposition}

The proof of this fact follows by noting that $S(z)$ is a
meromorphic doubly periodic function of $z$ with at most double poles
at the points $z_j$. By
expanding $S(z)$ in a Laurent series up to the constant term at
the $z=z_j$, we find that the difference
between left-hand side and right-hand side is a differential
operator whose coefficients are regular elliptic functions
vanishing  at least at one point. Such an operator vanishes
by Liouville's theorem.

The eigenvalue problem for common eigenfunctions of 
$H_0,\dots,H_n$ can then be formulated as $S(z)u=q(z)u$
with $q(z)=\sum_{k=1}^n\frac{c_k}2\bar\wp(z-z_k)
              +\sum_{k=1}^n\epsilon_k\bar\zeta(z-z_k)
              +\epsilon_0$. The eigenvalue of $H_j$ is then
$\epsilon_j$, $j=0,\dots,n$. Since $\sum_{j\geq 1} H_j=0$, one
must necessarily have $\epsilon_1+\cdots+\epsilon_n=0$.
        
Common eigenfunctions of $H_j$ and thus of $S(z)$ can
be obtained by the Bethe ansatz method. They have
the form $f(w_1)\cdots f(w_m)v_0$ where $v_0$ is the
tensor product of  highest weight vectors, $m=\frac12\sum\Lambda_j$
and $w_1,\dots,w_m$ are a solution to the system of Bethe
ansatz equations, see \cite{FV1,FV2} and below.

\subsection{Separation of variables}\label{ss-32}
We realize the representations of $sl_2(\C)$ by differential
operators:

\begin{lemma}\label{l-Verma} For any $\Lambda\in\C$,
the map $f\mapsto t$, $h\mapsto -2t\frac d{dt}+\Lambda$,
$e\mapsto -t\frac {d^2}{dt^2}+\Lambda\frac d{dt}$
defines a representation of $sl_2(\C)$ on $\C[t]$, the
Verma module $M_\Lambda$. If
$\Lambda$ is a nonnegative integer, $t^{\Lambda+1}\C[t]$
is an invariant subspace and the quotient $L_\Lambda=
\C[t]/t^{\Lambda+1}\C[t]$ is irreducible with highest 
weight vector $1\in \mathrm{Ker}(e)$ of weight $\Lambda$.
\end{lemma}

The proof consists of checking the relations
$[e,f]=h$, $[h,e]=2e$, $[h,f]=-2f$ and $e\, t^{\Lambda+1}=0$.

\medskip

Therefore we may realize the tensor product $M$
as $\C[t_1,\dots,t_n]$, and the tensor product of
irreducible representations (for integer $\Lambda_j$)
as $\C[t_1,\dots,t_n]/\sum_j (t_j^{\Lambda_j+1}\C[t_1,\dots,t_n])$.
Then $M[0]$ consists of homogeneous polynomials
of degree $m=\sum\Lambda_k/2$. We may then view
$e_\lambda(z),f_\lambda(z),h(z)$,
$ S(z)$, $H_j$ as  differential operators in $n+1$ variables
$\lambda,t_1,\dots,t_n$. They commute when acting on functions
which are homogeneous in $t_1,\dots, t_n$ of degree $m$.

\subsection{Separated variables}
We express the differential operators $S(z)$ in terms of
new variables so that the eigenvalue problem is reduced
to ordinary differential equations. 
Following Sklyanin's idea, the new variables $C,y_1,\dots,y_n$ are the
zeros and the leading coefficient of the operator $f_\lambda$:
\[
f_\lambda(z)=C\prod_{j=1}^n\frac{\theta(z-y_j)}{\theta(z-z_j)}
\]
Since $f_\lambda(z)$ is realized as a multiplication operator
and both sides of this equation are functions of $z$ 
with definite transformation properties under translations
by the lattice $\Gamma$, this equation does define, locally
around a generic point, a
biholomorphic change of variables $(C,y_1,\dots,y_n)\mapsto
(\lambda,t_1,\dots,t_n)$. The formulae are
\[
t_i=C\frac
{
\prod_{j}\theta(z_i-y_j)}
{
\theta'(0)\prod_{j:j\neq i}\theta(z_i-z_j)},
\]
and
\[
\lambda=\sum_{j=1}^n(y_j-z_j).
\]
{}From these formulae we deduce the transformation properties
of partial derivatives:
\[
\frac{\partial}
{\partial y_j}
=\frac{\partial}
{\partial\lambda}
+
\sum_{k=1}^n
\frac{\theta'(y_j-z_k)}
{\theta(y_j-z_k)}\, t_k\frac\partial{\partial t_k}\,.
\]
\[
C\frac\partial
{\partial C}=
\sum_{k=1}^n t_k
\frac\partial{\partial t_k}.
\]
The next step is to remark that a function $u(C,y_1,\dots,y_n)$
obeys $S(z)u=q(z)u$ with $q(z)=\sum_{k=1}^n\frac{c_k}2\bar\wp(z-z_k)
              +\sum_{k=1}^n\epsilon_k\bar\zeta(z-z_k)
              +\epsilon_0$ and $\sum_{j\geq 1}\epsilon_j=0$ if and
only if it obeys $S(y_j)u=q(y_j)u$ for all $j=1,\dots, n$ and all
generic points $(y_1,\dots,y_n)$. Here the ambiguous notation $S(y_j)$
means: write differential operator $S(z)$ with the coefficients on the
left of the partial derivatives and replace $z$ by $y_j$ in the
coefficients. To prove this statement notice that, if
$S(y_j)v=q(y_j)v$, then
$\prod\theta(z-z_i)(S(z)-q(z))u(y_1,\dots,y_n)$ 
 is a holomorphic theta function in $z$ of order $n$ vanishing
at $n$ generic points $y_i$. It thus vanishes, see the Appendix.

It is then convenient to use the identity
$
[e_\lambda(z),f_\lambda(z)]=-h'(z)$ on $M[0]$, to write $S(z)$ as
$
S(z)=\left(\frac{\partial}{\partial\lambda}-\frac12
h(z)\right)^2-h'(z)+2f_\lambda(z)e_\lambda(z).
$
so that the last term vanishes if we set $z=y_j$ and we get
\begin{equation}\label{e-sj}
S(y_j)=
\left(\frac{\partial}{\partial y_j}-\sum_{k=1}^n
\frac{\Lambda_k}2
\bar\zeta(y_j-z_k)
\right)^2.
\end{equation}
\begin{proposition}\label{p-88} A function $u(\lambda,t_1,\dots,t_n)$,
homogeneous of degree $m=\frac12\sum\Lambda_k$ in the
$t_i$,
is a local solution of the partial differential equations 
$S(z)u=q(z)u$, $z\in\C$ if and only if 
\[
u(\lambda,t_1,\dots,t_n)=C^nv(y_1,\dots,y_n)
\]
and $v$ obeys
\begin{equation}\label{sepeq}
\nabla_{y_j}^2v=q(y_j)v, \qquad \nabla_{y_j}=
\frac{\partial}{\partial y_j}-\sum_{k=1}^n
\frac{\Lambda_k}2
\bar\zeta(y_j-z_k).
\end{equation}

\end{proposition}

\subsection{Interpolation formula}
The formula \Ref{e-sj} expresses the values at $y_1,\dots,y_n$
of the coefficients
of $S(z)$ for $z=y_j$. Since the coefficients 
of $S(z)-\frac12\sum_{k}c_k\bar\wp(z-z_k)$ are elliptic
functions of $z$ with at most simple 
 poles at  $z_1,\dots,z_n$, they are uniquely determined
by these values, and can be calculated by an interpolation
formula: let us write $S(z)$ in the form
\[
S(z)=\frac12
\sum_{k=1}^nc_k
\bar\wp(z-z_k)+\prod_{k=1}^n\theta(z-z_k)^{-1}
\hat S(z),
\]
so that $\hat S(z)$ is a theta function of order $n$. Thus
(see the Appendix)
\[
\hat S(z)=\sum_{i=1}^n
\frac{\theta(z+\sum_{j\neq i}y_j-\sum_kz_k)}
{\theta(\sum_jy_j-\sum_kz_k)}\prod_{j:j\neq i}
\frac{\theta(z-y_j)}{\theta(y_i-y_j)}\hat S_i,
\]
with
\[
\hat S_i=\prod_{k=1}^n\theta(y_i-z_k)
\left(\left(\frac{\partial}{\partial y_i}-\sum_{k=1}^n
\frac{\Lambda_k}2
\bar \zeta(y_i-z_k)\right)^2
-\frac12\sum_{k=1}^nc_k\bar\wp(y_i-z_k)\right).
\]
\subsection{Bethe ansatz}\label{ss-BA}
The separated equation reads
\begin{equation}\label{sep}
\nabla_y^2v-(\sum_{k=1}^n\frac{c_k}2\bar\wp(y-z_k)
              +\sum_{k=1}^n\epsilon_j\bar\zeta(y-z_k)
              +\epsilon_0)v=0,
\end{equation}
with $\nabla_y=\frac{\partial}{\partial y}-\sum_{k=1}^n
\frac{\Lambda_k}2
\bar\zeta(y-z_k)
$.
It is a second order ordinary differential equation
with regular singular points at $z_k$ and characteristic
exponents $0$ and $\Lambda_k+1$. Following Hermite's
method to solve the Lam\'e equation (see \cite{WW}) we seek solutions of
the form
\begin{equation}\label{BA}
v(y)=e^{cy}\prod_{k=1}^m\theta(y-w_k).
\end{equation}
Functions of this form are called elliptic polynomials,
see the Appendix.

Let us first assume that $w_k\neq z_j \mod \Gamma$ for all $j,k$.
Then we also have $w_k\neq w_l$ for $k\neq l$, since the only
solution vanishing with its derivative at a regular point is
the trivial solution.
Rewrite the equation in the form
$v''(y)-\sum_k\Lambda_k\bar\zeta(y-z_k)v'(y)+b(y)v(y)=0$ so that
$b(y)$ has at most simple poles at the $z_j$.
Taking derivatives of a function $v$ of the form \Ref{BA}
and setting $y$ equal to one of its zeros, we find the
relation:
\[
v''(w_k)=\sum_{j\neq k}\bar\zeta(w_k-w_j)v'(w_k)+2cv'(w_k).
\]
Inserting this into the differential equation, we see that
$v$ is a solution if and only 
if its zeros $w_j$ obey the ``Bethe ansatz equations''
\[
\sum_{l=1}^n\Lambda_l\bar\zeta(w_j-z_l)-\sum_{k:k\neq j}\bar\zeta(w_j-w_k)=
2c,
\qquad j=1,\dots,m.
\]
Let us now consider the more general case of elliptic polynomials \Ref{BA}
vanishing at $z_i$, for $i$ in some subset $I$ of $\{1,\dots,n\}$.
Since the characteristic exponents at $z_i$ are
$0,\Lambda_i+1$, a solution vanishing at $z_i$ must vanish to order
$\Lambda_i+1$ and is thus divisible by $\theta(y-z_i)^{\Lambda_i+1}$
(this is only possible if $\Lambda_i\in\Z_{\geq0}$). 
Then $\tilde v(y)=\prod_{i\in I}\theta(y-z_i)^{-\Lambda_i-1}v(y)$ is
again of the form \Ref{BA}, but with $m$ replaced by $\tilde m
=m-\sum_{i\in I}(\Lambda_i+1)$. It obeys the equation \Ref{sep}
with $\Lambda_i$ replaced by $-\Lambda_i-2$ for $i\in I$.

Thus all elliptic polynomial solutions of \Ref{sep} are
of the form 
\[
v(y)=e^{cy}\prod_{i\in I}\theta(y-z_i)^{\tilde\Lambda_i+1}
\prod_{k=1}^{\tilde m}\theta(y-w_k),
\]
for some subset $I$ of $\{j\,|\,\Lambda_j\in\Z_{\geq0}\}$, such that 
$w_j\neq w_l\neq z_i$ $(j\neq l)$ and $w_1,\dots,w_{\tilde m},c$
obey the Bethe ansatz equations
\begin{equation}\label{BaE}
\sum_{l=1}^n\tilde
\Lambda_l\bar\zeta(w_j-z_l)-\sum_{k:k\neq j}\bar\zeta(w_j-w_k)=
2c,
\qquad j=1,\dots,\tilde m.
\end{equation}
Here $\tilde\Lambda_l=-\Lambda_l-2$ if $l\in I$ and 
$\tilde\Lambda_l=\Lambda_l$ otherwise.

To each such solution there corresponds a common eigenfunction $u$
which, expressed in the separated variables,
is $u=C^n\prod v(y_i)$. Up to a nonzero constant we get
\begin{equation}\label{BEV}
u(\lambda)=e^{c\lambda}f(w_1)\cdots f(w_{\tilde m})v_I.
\end{equation}
Here $v_I=\prod_{i\in I}(f^{(i)})^{\Lambda_i+1}v_0$ is
a product of singular vectors. Only eigenvectors corresponding
to $I=\mathrm{\o}$ have  non-trivial projections
 to eigenvectors with values in the
tensor product of irreducible representations.

Eigenvectors of the form \Ref{BEV}, such
that $w_1,\dots,w_{\tilde m},c$ are a solution of the Bethe ansatz equations
\Ref{BaE} with $w_j\neq w_k \mod\Gamma$, $(j\neq k)$ and
$w_j\neq z_i\mod\Gamma$ are called
{\em Bethe eigenvectors}.

\subsection{Completeness of Bethe eigenvectors}\label{ss-compl}

Let us consider the common eigenvalue problem
\begin{equation}\label{eigPro}
H_iu(\lambda)=\epsilon_i u(\lambda), \qquad i=0,\dots,n,
\end{equation}

A natural class of functions preserved by the operators $H_i$
is given by meromorphic sections of a flat line bundle on $E$.
Namely, let for a character $\chi:\Gamma\to\C^\times$,
$\mathcal{H}(\chi)$ be the space of meromorphic functions $\lambda\mapsto
u(\lambda)\in M[0]$ such that $u(\lambda+1)=\chi(1)u(\lambda)$
and
\[
u(\lambda+\tau)=\chi(\tau)e^{\pi i\sum_jz_jh^{(j)}}u(\lambda).
\]
It is easy to see that $e_\lambda(z)$, $f_\lambda(z)$ and
$\partial/\partial\lambda-h(z)/2$ preserve functions with these
transformation properties, so that $S(z)$ and $H_j$ preserve
$\mathcal{H}(\chi)$.

It is then natural to look for eigenfunctions 
(non-trivial solutions of the differential equations \Ref{eigPro})
in $\mathcal{H(\chi)}$.

Let $\Sigma(\chi)$ be the set of $(\epsilon_0,\dots,\epsilon_n)\in
\C^{n+1}$ such that there exists a non-trivial function $u\in \mathcal
{H}(\chi)$ with $H_ju=\epsilon_ju$, $j=0,\dots,n$. Let
$\sigma\in\Gamma^*$ be the character such that 
$\sigma(r+s\tau)=(-1)^{r+s}$.

\begin{thm} Let $\chi\in\Gamma^*$.
Then $(\epsilon_0,\dots,\epsilon_n)\in\Sigma(\chi)$ if and only if
$\sum_{k\geq1}\epsilon_k=0$ and the separated problem 
\[
\nabla_y^2v(y)-\frac14\sum_k\Lambda_k(\Lambda_k+2)\bar\wp(y-z_k)v(y)
=(\sum_k\bar\zeta(y-z_k)\epsilon_k+\epsilon_0)v(y),
\]
 admits a non-trivial elliptic polynomial
solution $v\in\Theta_m(\sigma^m\chi)$.
In this case there is a Bethe eigenvector
with these eigenvalues.
\end{thm}

\noindent{\it Proof:} A common eigenfunction, viewed as a
 polynomial in $t_i$ has the form
\[
u(\lambda,t_1,\dots,t_n)=
\sum_{m_1+\cdots+m_n=m}u_{m_1,\dots,m_n}(\lambda)\prod_it_i^{m_i}.
\]
Replacing $t_i,\lambda$ as functions of the new variables, we
obtain
\[
u(\lambda,t_1,\dots,t_n)=C^nv(y_1,\dots,y_n),
\]
where $v(y_1,\dots,y_n)$ is a linear combination of products of
theta functions in each of the $y_j$
with coefficients $u_{m_1,\dots,m_n}(\sum y_j-\sum z_k)$
and $v(y_1,\dots,y_n)$ obeys, in each variable, the separated
second order equation \Ref{sepeq}. A priori
the meromorphic function $v(y_1,\dots,y_n)$
may have poles on the hyperplane $\sum y_i=\sum z_i$
mod $\Gamma$. However this is impossible: consider $v$ as a function
of, say $y_1$ with the other variable fixed at some generic position. Then $v$, as
a function of $y_1$, being a solution of a linear second order equation
may only have singularities at the poles $z_j$ of the coefficients.
   Moreover, since
$u$ is in $H(\chi)$, the functions $y_i\mapsto
v(y_1,\dots,y_n)$ belong to $\Theta_m(\sigma^m\chi)$.
Thus the separated problem admits a non-trivial
solution in $\Theta_m(\sigma^m\chi)$.
As shown in the previous subsection, such a solution gives
rise to a Bethe eigenvector.
$\square$

\section{The difference case}
\subsection
{Representations of the elliptic quantum group $E_{\tau,\eta}(sl_2)$}

The difference  version of the differential operators $e_\lambda(z), f_\lambda(z), \partial_\lambda-h(z)/2$,
are operators obeying the relations of the elliptic quantum group $E_{\tau,\eta}(sl_2)$.
Let us recall the definitions \cite{FV3}: we fix two complex parameters
$\tau,\eta$, such that Im$(\tau)>0$.
The definition of $E_{\tau,\eta}(sl_2)$ is based on a dynamical $R$-matrix
$R(z,\lambda)$ which we now introduce.
Let
\be
\alpha(z,\lambda)=
\frac{\theta(\lambda+2\eta)\theta(z)}
{\theta(\lambda)\theta(z-2\eta)},\qquad
\beta(z,\lambda)=
-\frac{\theta(\lambda+z)\theta(2\eta)}
{\theta(\lambda)\theta(z-2\eta)},
\ee
Let $V$ be a two dimensional
complex vector space with basis $e[1],e[{-1}]$, and 
let $E_{ij}e[k]=\delta_{jk}e[i]$, $h=E_{11}-E_{-1,-1}$.
 Then, for $z,\lambda\in\C$, $R(z,\lambda)\in\mathrm{End}(V\otimes V)$ 
is the matrix
\bea\label{exR}
 R(z,\lambda)&\!=\!&E_{11}\otimes E_{11}+E_{-1,-1}\otimes E_{-1,-1}+
 \alpha(z,\lambda)E_{11}\otimes E_{-1,-1}
\\
 &
+\!&\alpha(z,-\lambda)E_{-1,-1}\otimes E_{11}
+\beta(z,\lambda)E_{1,-1}\otimes E_{-1,1}
+\beta(z,-\lambda)E_{-1,1}\otimes E_{1,-1}.
\eea
This $R$-matrix obeys the dynamical 
 quantum Yang--Baxter equation
\bea
R^{(12)}(z\!-\!w,\lambda\!-\!2\eta h^{(3)})\!\!\!\!\!
&R^{(13)}(z,\lambda)\,
R^{(23)}(w,\lambda\!-\!2\eta h^{(1)})&\!\!\!\!\!
=
\\
  &R^{(23)}(w,\lambda)\,
R^{(13)}(z,\lambda\!-\!2\eta h^{(2)})&\!\!\!\!\!
R^{(12)}(z\!-\!w,\lambda)
\eea
in $\mathrm{End}(V\otimes V\otimes V)$, $z,w,\lambda\in\C$.  The meaning of this
notation is the following: $R^{(12)}(\lambda-2\eta h^{(3)})\,v_1\otimes
v_2\otimes v_3$
is defined as
\be
(R(z,\lambda-2\eta\mu_3)\,v_1\otimes v_2)\otimes v_3,
\ee
if $h v_3=\mu_3v_3$.
The other terms are defined similarly: in general, let 
$V_1,\dots, V_n$ be modules over the one dimensional Lie
algebra $\h=\C h$ with one generator $h$, such that, for all $i$,
$V_i$ is the direct sum of finite dimensional eigenspaces
$V_i[\mu]$  of $h$, labeled by the eigenvalue $\mu$.
We call such modules diagonalizable $\h$-modules.
If $X\in\mathrm{End}(V_i)$
we denote by $X^{(i)}\in\mathrm{End}(V_1\otimes\cdots\otimes V_n)$ 
the operator $\cdots\otimes\Id\otimes X\otimes\Id\otimes\cdots$
acting non-trivially on the $i$th factor, and
if $X=\sum X_k\otimes Y_k\in\mathrm{End}(V_i\otimes V_j)$ we set
$X^{(ij)}=\sum X_k^{(i)}Y_k^{(j)}$. If $X(\mu_1,\dots,\mu_n)$
is a function with values in $\mathrm{End}(V_1\otimes\cdots\otimes V_n)$,
then $X(h^{(1)},\dots,h^{(n)})v=X(\mu_1,\dots,\mu_n)v$ if
$h^{(i)}v=\mu_iv$, for all $i=1,\dots,n$.

 A {\em module over} $E_{\tau,\eta}(sl_2)$ is then
a diagonalizable $\h$-module $W=\oplus_{\mu\in\C}  W[\mu]$,
together with an $L$-operator $L(z,\lambda)\in\mathrm{End}_\h(V\otimes W))$
(a linear map commuting with $h^{(1)}+h^{(2)}$)
depending meromorphically on  $z,\lambda\in\C$ and obeying the
relations
\begin{eqnarray}\label{e-RLL}
R^{(12)}(z\!-\!w,\lambda\!-\!2\eta h^{(3)})\!\!\!\!\!
&L^{(13)}(z,\lambda)\,
L^{(23)}(w,\lambda-2\eta h^{(1)})&\!\!\!\!\!
=
\\
  &L^{(23)}(w,\lambda)\,
L^{(13)}(z,\lambda\!-\!2\eta h^{(2)})&\!\!\!\!\!
R^{(12)}(z\!-\!w,\lambda)\,.\nonumber
\end{eqnarray}
For example, $W=V$, $L(w,\lambda)=R(w-z_0,\lambda)$ 
is a module over $E_{\tau,\eta}(sl_2)$,
called the fundamental representation, with evaluation point $z_0$.
In \cite{FV3} more general examples of such modules were constructed: in 
particular, for any pair of complex numbers $\Lambda, z$ we
have an {\em evaluation Verma module} $M_\Lambda(z)$. It has
a weight decomposition $M_\Lambda=\oplus_{j=0}^\infty M_\Lambda[\Lambda-2j]$,
with one-dimensional weight spaces $M_{\Lambda}[\mu]$. The action
of the $L$-operator
is described explicitly in \cite{FV3}.
Also, we have a notion of tensor products of modules over
$E_{\tau,\eta}(sl_2)$. The main examples considered in this
paper will be tensor products $M_{\Lambda_1}(z_1)\otimes
\cdots \otimes M_{\Lambda_n}(z_n)$ of evaluation Verma modules and
some of their subquotients.

It will be convenient here to consider more general $L$-operators obeying
the relations. So we define a {\em functional module} over $E_{\tau,\eta}(sl_2)$
to be given by a pair $(W,L)$ where $W$ is a space of complex-valued functions
on a certain set and $h$ acts on it as multiplication by a function, and
$L(z,\lambda)$ is a meromorphic function of $z$ and $\lambda$ acting as
a difference operator on $V\otimes W$, commuting with $h\otimes1+1\otimes h$,
and obeying the relations \Ref{e-RLL}.
An example of such functional modules is provided by the
``universal evaluation modules'' of \cite{FV3}. $h$ acts by multiplication
by a continuous variable. Evaluation Verma modules are obtained by 
restricting the range of this continuous variable to a discrete set.

For any module or functional module 
$W$ over $E_{\tau,\eta}(sl_2)$,  we define the
associated operator algebra, an algebra of operators on the space $\mathrm{Fun}(W)$
of meromorphic functions of $\lambda\in\C$
 with values in $W$. It is  generated by $h$,
acting on the values, and
operators  $a(z), b(z), c(z), d(z)$. Namely, let $\tilde L(z)
\in \mathrm{End}(V\otimes\mathrm{Fun}(W))$ be the operator
defined by 
$(\tilde L(z)(v\otimes f))(\lambda)=
L(z,\lambda)(v\otimes f(\lambda-2\eta\mu))$ 
if $h v=\mu v$. View $\tilde L(z)$
as a 2 by 2 matrix  with entries in $\mathrm{End}(\mathrm{Fun}(W))$:
\bea
\tilde L(z) (e[1]\otimes f)&\!=\!&e[1]\otimes a(z)f+e[{-1}]\otimes c(z)f,
\\
\tilde L(z) (e[{-1}]\otimes f)&\!=\!&e[1]\otimes b(z)f+e[{-1}]\otimes d(z)f.
\eea
The relations obeyed by these operators are described in detail
in \cite{FV3} (in \cite{FV3} these operators are denoted by 
$\tilde a(z), \tilde b(z)$ and so on).

To each module we associate a central element of the operator
algebra. It is given by the quantum determinant \cite{FV3}
\[
\mathrm{Det}(z)
=
\frac
{\theta(\lambda)}
{\theta(\lambda-2\eta h)}
(a(z+2\eta)d(z)-c(z+2\eta)b(z)).
\]

\subsection{A class of representations by difference operators}
\label{s-diffop}

 Let $z_1, \ldots , z_n \in \C$ be distinct points
     and $\Lambda_1, \ldots, \Lambda_n \in Z_{\geq0}$.
Let us introduce difference  operators acting on functions of $n+1$ complex variables
$\lambda,x_1,\dots,x_n$. Let $(T_{x_i}^{a}f)(\lambda,x_1,\dots,x_n)=f(\lambda,x_1,\dots,x_i+a,\dots,x_n)$
and $(T_\lambda^af)(\lambda,x_1,\dots,x_n)=f(\lambda+a,x_1,\dots,x_n)$. The steps $a$ will always
be $\pm2\eta$. Let 
\[\Delta_{-} (z) = \prod_{i=1}^{n}
\theta(z - z_i + \Lambda_i \eta) \mbox{ and } \Delta_{+}(z) = \prod_{i=1}^{n} 
\theta(z - z_i - \Lambda_i \eta) \, .
\]
The functions $\Delta_{\pm} (z=-x_i)$,
considered as multiplication operators,
 will be denoted simply by $\Delta_{\pm}(-x_i)$. We also set 
$s=\sum_{i=1}^n(x_i+z_i)$.
 
With these conventions, we define: 
                     \begin{eqnarray*}
                      & a(z) & = \prod_{i=1}^{n} \theta(z+x_i) \frac{ 
                      \theta(\lambda + \sum_{l=1}^{n}(x_l + z_l + \Lambda_l 
                      \eta)) }{ \theta(\lambda)} T_{\lambda}^{- 2 \eta}\\ 
                      & b(z) & = -\sum_{i=1}^{n} \frac{\theta(
                      \lambda + z + x_i)}{\theta(\lambda)} \prod_{j \neq i}
                     \frac{\theta(z+x_j)}{\theta(x_i - x_j)} \Delta_{+}(-x_i)
                       T_{x_i}^{- 2 \eta} T_{\lambda}^{+2 \eta}, \\
                      & c(z) & =  -\sum_{i=1}^{n} \frac{\theta(
                      - \lambda + z + x_i - 2 s 
                       )}{\theta(\lambda)} \prod_{j \neq i}
                      \frac{\theta(z+x_j)}{\theta(x_i - x_j)} \Delta_{-}(-x_i)
                       T_{x_i}^{+2 \eta} T_{\lambda}^{-2 \eta}, \\
                      & \mathrm{Det}(z) & = \prod_{i=1}^{n}
                      \theta(z - z_i-\Lambda_i \eta) \, \theta(z - z_i+ 
                      \Lambda_i \eta +2 \eta)    
                     \end{eqnarray*} 
\begin{thm}\label{p-0.0}
The difference operators $a(z),b(z),c(z)$ together with $d(z)$ defined implicitly
by the determinant relation
\[
a(z+2\eta)d(z)-c(z+2\eta)b(z)=\frac{\theta(\lambda-2\eta h)}{\theta(\lambda)}\mathrm{Det}(z)
\]
obey the relations of the elliptic quantum group \Ref{e-RLL} with
$\eta h=-\sum_i(x_i+z_i)$.
\end{thm}

\noindent{\bf Example.} If $n=1$, the generators act on functions
of two variables $x_1,\lambda$, and $h$ acts as $-x_1-z_1$. If we introduce
a new variable $h=-\eta^{-1}(x_1+z_1)$, so that the generator $h$ acts
by multiplication by $h$, we obtain a representation on functions of $h,\lambda$.
the action of the generators is given by the difference operators
       \begin{eqnarray*}
                       a(z)v(h,\lambda) & = &\theta(z-z_1-\eta h) \frac{ 
                      \theta(\lambda -\eta h + \Lambda_1 
                      \eta) }{ \theta(\lambda)} 
       v(h,\lambda-2\eta)\\ 
                      b(z) v(h,\lambda) & = & \frac{\theta(
                      \lambda + z -z_1-\eta h)}{\theta(\lambda)} 
                      \theta(-\eta h + \Lambda_1\eta)
       v(h+2,\lambda+2\eta),\\
                      c(z) v(h,\lambda) & = & -  \frac{\theta(
                      - \lambda + z  -  z_1 +\eta h
                       )}{\theta(\lambda)} 
                   \theta(\eta h+ \Lambda_1 \eta)
       v(h-2,\lambda-2\eta),\\
                       d(z) v(h,\lambda) & = & 
{\theta (z-z_1+ \eta h)} 
\frac{\theta(\lambda - \eta h  - 
                    \Lambda_1 \eta)}{ \theta(\lambda)}
       v(h,\lambda+2\eta).
\end{eqnarray*} 
This is, up to normalization, the ``universal evaluation module''
of \cite{FV3}, Sect.\ 9.
\medskip

\noindent{\it Proof of Theorem \ref{p-0.0}:}
The proof consists of a straightforward verification of the sixteen relations. Let
us give an example: one relation is
         \begin{eqnarray*}
         a(z)b(w) = \frac{\theta(z-w)\theta(\lambda+ 2\eta)}{\theta(z-w-2\eta)
                    \theta(\lambda)} b(w) a(z) + \frac{\theta(z-w-\lambda)   
                    \theta(2\eta)}{ \theta(z-w-2\eta)\theta(\lambda)}a(w)b(z)
                    \,. 
        \end{eqnarray*}
         This identity is verified
 by looking at each summand of $b(z)$ and $b(w)$ 
         separately. The corresponding equation to one such summand typically 
         looks like
        \begin{eqnarray*}
        &  & \theta(z+x_k) (\prod_{i \neq k} \theta(z+x_i)) \frac{\theta( 
        \lambda + 
        \sum_{i=1}^n(x_i + z_i +\Lambda_i \eta)) \theta(\lambda + w + x_k - 2 
        \eta) }{\theta(\lambda) \theta(\lambda - 2\eta)}\times \\
        & & \times ( \prod_{i \neq k}
        \frac{\theta(w+x_i)}{\theta(x_k - x_i)}) \Delta_{+}(-x_k)
         T_{x_k}^{-2\eta}\\
        & = & \frac{\theta(z-w)\theta(\lambda+ 2\eta)}{\theta(z-w-2\eta)
        \theta(\lambda)} \frac{\theta(\lambda + w + x_k)}{\theta(\lambda)}
        (\prod_{i \neq k} \frac{\theta(w+x_i)}{\theta(x_k - x_i)}) \times \\
        & & \times
        (\prod_{i \neq k} \theta(z+x_i)) \frac{\theta(z+x_k - 2\eta) \theta( 
        \lambda + \sum_{i=1}^n(x_i + z_i +\Lambda_i \eta))}{\theta(\lambda + 
       2\eta)} \Delta_{+}(-x_k) T_{x_k}^{-2\eta} \\
       & + & \frac{\theta(z-w-\lambda) \theta(2\eta)}{ \theta(z-w-2\eta) 
        \theta( \lambda)} \frac{\theta(\lambda + \sum_{i=1}^n(x_i + z_i + 
        \Lambda_i \eta)) \theta(w+x_k)}{\theta(\lambda)}(\prod_{i \neq k} 
        \theta(w+x_i)) \times \\
       & & \times (\prod_{i \neq k} \frac{\theta(z+x_i)}{\theta(x_k - x_i)})
        \frac{\theta(z+x_k+\lambda-2\eta)}{\theta(\lambda - 2 \eta)} 
        \Delta_{+}(-x_k) T_{x_k}^{-2\eta} \, .  
        \end{eqnarray*}   
       By taking into account that each summand of the above equation
       involves a factor
       \begin{eqnarray*}
       (\prod_{i \neq k}\frac{\theta(z+x_i)\theta(w+x_i)}{\theta(x_k - x_i)})
       \frac{\theta(\lambda + \sum_{i=1}^n(x_i + z_i + 
        \Lambda_i \eta))}{\theta(\lambda)} \, ,
       \end{eqnarray*}  
       the task reduces to verifying
       \begin{eqnarray*}
       & & \theta(z+x_k) \frac{\theta(\lambda+w+x_k-2\eta)}{\theta(\lambda - 
       2\eta)} = \frac{\theta(z-w)\theta(\lambda + w + x_k) \theta(z+x_k - 
       2\eta)}{\theta(z-w-2\eta) \theta(\lambda)}  \\
       & & + \frac{\theta(z-w-\lambda)
       \theta(2 \eta) \theta(w+x_k) \theta(\lambda + z + x_k - 2\eta)}{ \theta        (z- w- 2\eta) \theta( \lambda) \theta(\lambda -  2\eta)} \, , \\
       & & \mbox{ which we write in the form } \\
       & & f_1 (z,w,\lambda) = f_2 (z,w,\lambda) + f_3 (z,w,\lambda) \, .
       \end{eqnarray*} 
       This identity is proved in two steps.
       First, one shows that the functions $f_i (z,w,\lambda),$ $ i=1,2,3 \, ,$
       transform in the same way under $\lambda \rightarrow \lambda+ 1, 
\lambda \rightarrow
       \lambda  + \tau \ , $. The transformation
       laws thus obtained are the following:
       \begin{eqnarray*} f_i (z,w,\lambda+1) = f_i (z,w,\lambda) \\
                         f_i (z,w,\lambda+\tau) = 
e^{-2\pi i(w+x_k)} f_i (z,w,\lambda).
                         \end{eqnarray*}
       for $i=1,2,3 \,.$ 
       The second step is
 to show that the above relation holds for the residues of
       the functions $f_i(z,w,\lambda)$. Here,
       one has to show that the identity
       holds for $\lambda = 2 \eta$ and $\lambda = 0 \, .$
       For $\lambda=0$ one obtains
       \begin{eqnarray*}
       \frac{\theta(z-w) \theta(w+x_k) \theta(z+x_k-2\eta)}{\theta(z-w-2\eta)}
       + \frac{\theta(z-w) \theta(2\eta) \theta(w+x_k) \theta(z+x_k - 2\eta)}{ 
       \theta(- 2\eta) \theta(z-w-2\eta)} = 0 \, , 
       \end{eqnarray*}
       whereas the $\lambda = 2\eta$ residue yields
       \begin{eqnarray*}
       \theta(z+x_k) \theta(w+x_k) = \frac{\theta(z-w-2\eta) \theta(2 \eta)}{ 
       \theta(z-w-2\eta) \theta(2 \eta)} \theta(w+x_k) \theta(z+x_k) \, .
       \end{eqnarray*}
       This proves ~\cite{FV3}, paragraph 3, relation 2. The other relations
       are proved in a similar but often more 
       intricate fashion.

One identity that is used is the vanishing of the 
sum of the residues at $v=-x_j-2\eta$, $-x_j$ of the
function
\[
f(v)=
\frac {\theta(2s + x_i + v + 2 \eta)}{ 
           \theta ( v + x_i + 2 \eta)} \prod_{l=1}^n \frac{ \theta (v - z_l - 
           \Lambda_l) \theta (v - z_l + \Lambda_l + 2\eta)}{\theta (v + x_l) 
           \theta(v + x_l + 2\eta)} \, , 
\]
a consequence of the double periodicity $f(v+1)=f(v+\tau)=f(v)$.

Also some of the more complicated relations, such as $d(z)d(w)=d(w)d(z)$
turn out to be consequences of the simpler relations and the fact
that the determinant is central.
$\square$

\medskip

\noindent{\bf Remark.} This is the difference elliptic analogue of formulae
that have appeared in the literature. In the rational difference and differential case such
a formula has been written by Sklyanin \cite{S}. A trigonometric
difference version appears in Tarasov and Varchenko \cite{TV}. 

\subsection{Restrictions}
The linear difference operators defined in the previous subsections
act on meromorphic funtions of complex variables $\lambda, x_1,\dots, x_n$. To compare
these operators with evaluation modules of the elliptic quantum
groups and with transfer matrices of IRF models, we have to restrict their
action 
to the space of
meromorphic function defined on submanifolds of $\C^{n+1}$. The conditions
for a difference operator $X$ with meromorphic coefficients
to be defined on meromorphic functions on a submanifold $S$ are
that the value at a generic $x\in S$ of $Xf(x)$ is well-defined
(i.e., there are no poles at generic points of $S$) and
is only a function  of the values of $f$ at
points of $S$. Equivalenty,  a difference operator
$X$ can be restricted to $S$ if it 
maps functions vanishing on $S$ to functions vanishing
on $S$. The restriction is then identified with the induced
action on the quotient by the function vanishing on $S$.

These conditions are fulfilled in the following situations:
\begin{enumerate}
\item{\em Restriction to discrete values of $x_i$.} We assume
that $z_i$, $\eta$ are generic and that the $\Lambda_i$ are
non-negative integers. We take $S$ to be the set
\begin{eqnarray*}
S_0&=&\{(\lambda,x_1,\dots,x_n)\in\C^{n+1}\,
|\, -x_i=z_i+\eta(\Lambda_i-2m_i), \\
&&m_i=0,1,\dots,\Lambda_i, 
i=1,\dots,n\}.\nonumber
\end{eqnarray*}
Since the steps in the difference operators are by multiples
of $2\eta$, it is clear that one can restrict the action
of $a(z),\dots,d(z)$ to functions on subsets
given by conditions $x_i\in a_i+2\eta\Z$, for generic $a_i,\eta$.
The genericity condition on $a_i$, $\eta$ 
is that the poles at $x_i-x_j$ mod
$\Gamma$ of the coefficients of the difference operators are
never on the subset. What we have 
to check is that the restriction to these finite
sets of values for $x_i$ is well-defined. Since only $b(z)$,
$c(z)$ shift the value of $x_i$, it is sufficient to consider
these two operators.  For a function $f(\lambda,x_1,\dots,x_n)$,
the value at $-x_i=z_i+\eta\Lambda_i$ of $b(z)f$ appears to 
depend on the value of 
$f$ at $-x_i+2\eta=z_i+\eta\Lambda_i+2\eta$, which is not in $S_0$, 
but in fact
it does not, since the coefficient $\Delta_-(-x_i)$ vanishes there.
Similarly $c(z)$ is well-defined on $S_0$.

For this restriction $h$ has discrete spectrum. It will be
used to compare our representation with tensor products
of irreducible representations.

\item{\em Restriction to $\lambda=\eta h$.}
Let $S$ be the set
\[
S_1=\{(\lambda,x_1,\dots,x_n)\,|\, \lambda=-\sum_i(x_i+z_i)\}.
\]
Then $b(z),c(z)$ can be restricted to functions on $S_1$. Indeed,
if $u$ is a function vanishing on $S_1$, then $T_{x_i}^{-2\eta}
T_\lambda^{2\eta}u$ still vanishes on $S_1$. The denominators
of the coefficients of $b(z)$ do not vanish at generic points
of $S_1$. Thus $b(z)u=0$ if $u$ vanishes on $S_1$.  
The same reasoning applies
to $c(z)$.

This restriction is needed, as we shall see, to construct commuting
transfer matrices.

\item{\em IRF restriction.} Consider the restriction
of $b(z), c(z)$ on functions on $S=S_0\cap S_1$. If $z_i$ and $\eta$
are generic, the only possible pole on $S$ in the coefficients
of these differential operators comes from the 
denominator $\theta(\lambda)$. This denominator does
not vanish if we assume for instance that the $\Lambda_i$
are all odd.

Here $S$ is finite, so that the restriction is to a finite
dimensional space of functions, which will be identified with
the space of states of an IRF model.
\end{enumerate}

\subsection{Commuting difference operators}\label{s-twisted}
One of the main properties of $R$-matrices to statistical
mechanics is that
$L$-operators obeying quantum group relations give rise to
commuting transfer matrices $\mathrm{tr}_VL$. In \cite{S} it
is noticed that more generally one can consider $\mathrm{tr}_V((K
\otimes 1)L)$
for some endomorphism of $V$ such that $K\otimes K$ commutes with
the $R$-matrices. 

In our dynamical
case it is known that the traces $a(z)+d(z)$ commute for different
values of $z$ when acting on the zero weight space of a module
over $E_{\tau,\eta}(sl_2)$. Another possibility to obtain commuting
operators is to take the {\em twisted transfer matrix}
$\mathrm{tr}_V((K\otimes 1)L)$ with a suitable $K$:

\begin{proposition}\label{p-twisted} For any $\vartheta\in\C$,
the operators $T(z)=b(z)+\vartheta c(z)$, $z\in\C$,
restricted to functions on the submanifold $S_1$
given by the equation $\lambda+\sum_{i=1}^nx_i+\sum_{i=1}^nz_i=0$,
form a commuting family.
\end{proposition}

This proposition can be proved directly or by the following general
argument.

It is first of all sufficient to consider the case $\vartheta=1$
since the other cases are obtained by conjugating the operator by
the multiplication by an exponential function of $\lambda$.
Then we may write 
$T(z)$ as
\[
T(z)=\sum_{\mu=\pm1}
\tr_{V[\mu]}(K\otimes 1\;L(z,\lambda))T_\lambda^{-2\eta\mu}\,,
\]
with
\[
 K  =  \left( \begin{array}{cc} 0 & 1 \\ 1 & 0 
   \end{array} \right) \, .
\]
The partial trace $\tr_{V[\mu]}:\mathrm{End}(V\otimes W)\to 
\mathrm{End}(W)$ is the homomorphism 
such that if $a\in
\mathrm{End}(V)$, $b\in
\mathrm{End}(W)$,
then $\tr_{V[\mu]}(a\otimes b)=\sum e^i(a e_i)b$ for any basis
 $(e_i)$ of $V[\mu]$ with $e^i\in V^*$ defined by $e^i(e_j)=\delta_{ij}$
and $e^i(w)=0$ for $w\in V[\nu]$, $\nu\neq\mu$.
We have
\begin{equation}\label{e-98}
K\otimes K\, R(z,\lambda)=R(z,-\lambda)\,K\otimes K.
\end{equation}
We then write the RLL relations in the form
\begin{eqnarray*}
 &  R^{(12)} (z-w,\lambda - 2 \eta h) 
   L^{(1)} (z,\lambda) L^{(2)} (w,\lambda - 2 \eta h^{(1)})
   R^{(12)} (z-w,\lambda)^{-1}  \, = & \\
   & L^{(2)} (w, \lambda) L^{(1)} (z, \lambda - 2\eta h^{(2)})
   \, , & 
\end{eqnarray*} 
multiply both sides by $(K\otimes K)^{(12)}$ from the
left, and take the partial trace over
a weight space $(V\otimes V)[\mu]$.
Using \Ref{e-98}, we obtain
\begin{eqnarray*}
 & \mathrm{tr}_{(V\otimes V)[\mu]}
( R^{(12)} (z-w,-\lambda +
   2 \eta h )K^{(1)}K^{(2)}  & \\
& L^{(1)}
   (z,\lambda) L^{(2)} (w,\lambda - 2 \eta h^{(1)}) R^{(12)} (z-w,\lambda)^{-1})
 \, = & \\
 & \tr_{(V\otimes V)[\mu]} 
K^{(2)} L^{(2)} (w,\lambda) K^{(1)}L^{(1)} (z,\lambda - 2 \eta 
    h^{(2)}). & 
\end{eqnarray*} 
The next step is to use the cyclicity of the trace to bring the first $R$ matrix
to the right. For this we need the commutation relations of $h$ with the product
in the trace. Since $[L^{(i)},h^{(i)}+h]=0$ and $h^{(i)}K^{(i)}=
-K^{(i)}h^{(i)}$,
we see that 
\[
h\, \tr_{(V\otimes V)[\mu]}(A^{(12)}K^{(1)}K^{(2)}L^{(1)}L^{(2)})=
\tr_{(V\otimes V)[\mu]}(A^{(12)}
K^{(1)}K^{(2)}L^{(1)}L^{(2)})\cdot (h+2\mu),
\]
for any $A\in
\mathrm{End}(V\otimes V)$ commuting with $h^{(1)}+h^{(2)}$.
We then get
\begin{eqnarray*}
 & \sum_\mu\mathrm{tr}_{(V\otimes V)[\mu]}
(K^{(1)}K^{(2)} L^{(1)}
   (z,\lambda) L^{(2)} (w,\lambda - 2 \eta h^{(1)})  & \\
& R^{(12)} (z-w,\lambda)^{-1}
 R^{(12)} (z-w,-\lambda +
   2 \eta (h+2\mu) )T_\lambda^{-2\eta\mu}
 \, = & \\
 & T(w)T(z) &
\end{eqnarray*}
If we then have a relation $\eta h\, u(\lambda)=\lambda u(\lambda)$, and we
apply the above equation to $u$, we may replace $h$ in the left-hand
side by $\lambda/\eta-2\mu$, and the $R$ matrices cancel, so that
$T(z)T(w)=T(w)T(z)$.

\subsection{Evaluation modules}
Here we show that the restriction of the operators $a,b,c,d$ of
Prop.\ \ref{p-0.0} to functions  on the submanifold
$S_0$
is essentially the tensor product of
finite dimensional irreducible evaluation modules of \cite{FV3}.

\begin{proposition}
Suppose that $\eta, z_1,\dots,z_n$ are generic.
Let $a,b,c,d$ be the difference operators defined in \ref{s-diffop} 
restricted
to functions on $S_{0}$ and let $\kappa(z)=\prod_{i=1}^n
\theta(z-z_i-\eta\Lambda_i)^{-1}$. Then
$\bar a(z)=\kappa(z)a(z),\dots,\bar d(z)=\kappa(z)d(z)$ define
an $E_{\tau,\eta}(sl_2)$ module isomorphic to the
tensor product $L_{\Lambda_1}(z_1-\eta)
\otimes\cdots\otimes L_{\Lambda_n}(z_n-\eta)$ of irreducible evaluation
modules.
\end{proposition}

Let $W$ be the space of functions on $S$. It is a
vector space over the field of meromorphic functions of 
$\lambda$ of dimension $\prod_{i=1}^n(\Lambda_i+1)$.
To prove this proposition we have to identify a highest weight
vector in $v$, i.e., an eigenvector of $a(z)$, $d(z)$ and
$h$ killed
by $c(z)$. The eigenvalues $(A(z,\lambda),D(z,\lambda),\Lambda)$
of $(a(z),d(z),h)$ determine then by \cite{FV3} uniquely an
irreducible module up to isomorphism.

Let $\delta_a(x_i)\in W$ be the delta function at $x_i=a$:
$\delta_a(x_i)=1$ if $x_i=a$, $\delta_a(x_i)=0$, if $x_i\neq a$.
     The highest weight vector may be taken as the product
of     delta functions
     $$ v_{h.w.} \, = \,
       \prod_{i=1}^n \delta_{-z_i - \Lambda_i \eta} (x_i) \, . $$  
This function is indeed annhilated by $c(z)$ and 
$hv_{h.w.}=\sum\Lambda_iv_{h.w.}$. Moreover $v_{h.w.}$ is an eigenvector
for $a(z)$ and $d(z)$:
     \begin{eqnarray*} 
     a(z)v_{h.w.}&=&
 \prod_{i=1}^{n} \theta ( z - z_i -
        \Lambda_i \eta) \frac{\theta(\lambda - \sum_{i=1}^{n} \Lambda_i \eta
       + \sum_{i=1}^{n} \Lambda_i \eta)}{\theta(\lambda)} v_{h.w.}
\\
& =&     \prod_{i=1}^{n} \theta ( z - z_i - \Lambda_i \eta) v_{h.w.}\, . 
     \end{eqnarray*} 
Thus $a(z)v_{h.w.}=A(z,\lambda)v_{h.w.}$, with $A(z,\lambda)=
     \prod_{i=1}^{n} \theta ( z - z_i - \Lambda_i \eta)$.
Similary $d(z)v_{h.w.}=D(z,\lambda)v_{h.w.}$ with eigenvalue
     \begin{eqnarray*}
 D(z,\lambda) &=&  \frac{\theta(\lambda -  2 \eta 
      \sum_{i=1}^n \Lambda_i )}{\theta ( \lambda)} Det(z - 2\eta ,\lambda) 
      A^{-1} (z-2\eta, \lambda+2\eta) \, . \\
&=& \frac{ \theta (
        \lambda - 2 \eta \sum_{i=1}^{n} \Lambda_i)}{\theta (\lambda)} 
         \prod_{i=1}^{n} \frac{\theta ( z - z_i - \Lambda_i \eta - 2 \eta)
        \theta ( z-z_i +\Lambda_i \eta )}{\theta 
        ( z - z_i - \Lambda_i \eta - 2 \eta)} \, . 
  \end{eqnarray*}  
  Thus the eigenvalues of $(\bar a(z),\bar d(z),h)$ are 
$(1,\bar D(z,\lambda),\sum\Lambda_i)$ with
     $$ \bar D(z, \lambda) = \frac{ \theta ( \lambda -
      2\eta
 \sum_{i=1}^n \Lambda_i)}{\theta (\lambda) } \prod_{i=1}^n \frac{ \theta
       (z-z_i + \Lambda_i \eta)}{\theta ( z-z_i - \Lambda_i \eta)} \, , $$
    which indeed reproduces the highest weight defined in ~\cite{FV3}, 
p.~750. 

\subsection{Separation of variables: continuous case}
In this section we find an analogue of the results of
\ref{ss-BA}, \ref{ss-compl}. We consider the continuous case,
 in which the variables $x_i$ take arbitrary complex values.
Then the eigenvalue problem for a function $u$ on $S_1$ reads
$T(z)u(x)=\epsilon(z)u(x)$, with
\begin{eqnarray*}
T(-z)u(x_1,\dots,x_n)&=&\sum_{i=1}^n
\frac{\theta(\lambda-z+x_i)}
{\theta(\lambda)}
\prod_{j\neq i}
\frac{\theta(z-x_j)}
     {\theta(x_i-x_j)}\nonumber\\
&&\bigl(\prod_{k=1}^n\theta(x_i+z_k+\eta\Lambda_k)
u(x_1,\dots,x_i-2\eta,\cdots,x_n)
\\
&&+\prod_{k=1}^n\theta(x_i+z_k-\eta\Lambda_k)
u(x_1,\dots,x_i+2\eta,\dots,x_n)\bigr)\,,\nonumber
\end{eqnarray*}
where we view a function $u$ on $S_1$ as a function of
$x_1,\dots,x_n$ by setting $\lambda=-\sum(x_j+z_j)$.

From this formula it is clear that $T(z)u(x)$ is an entire
holomorphic function of $z$ with theta function behavior as
$z$ is translated by elements of the lattice $\Gamma=\Z+\tau\Z$.
It follows that a necessary condition of $\epsilon(z)$ to be
an eigenvalue is that $\epsilon$ belong to the space 
$Theta_n(\chi_0)$ of theta functions of order $n$ with 
character $\chi_0:\Gamma\to\C^*$ (see the Appendix) given
by $\chi_0(r+s\tau)=(-1)^{n(r+s)}\exp(2\pi i\sum z_k)$. This
means that $\epsilon$ is an entire function obeying
\[
\epsilon(z+1)=\chi_0(1)\epsilon(z),\qquad
\epsilon(z+\tau)=\chi_0(\tau)e^{-\pi in(2z+\tau)}\epsilon(z).
\]
The method of separation of variables consists in looking
 for common eigenfunctions
$u(x)$ in the factorized
form $u(x)=\prod_{i=1}^n Q(x_i)$. Setting $z=-x_i$ in the
eigenvalue problem $(T(z)-\epsilon(z))\,u=0$ we see
that a necessary condition is that $Q,\epsilon$ obey
the difference equation
\[
A_+(x)\,
Q(x-2\eta)
+
A_-(x)\,
Q(x+2\eta)
=\epsilon(-x)Q(x), \qquad
A_{\pm}(x)=\prod_{k=1}^n\theta(x+z_k\pm\eta\Lambda_k)\,
\]
As explained in the Appendix, this difference equation
has an {\em elliptic polynomial solutions}, i.e., a solution
of the form 
\begin{equation}\label{e-BA2}
Q(x)=e^{ax}\prod_{k=1}^m\theta(x-w_k).
\end{equation}
if $\sum\Lambda_i$ is an even integer $2m$.
Such a solution may be constructed by the Bethe ansatz:

\begin{proposition}\label{p-ba}
Suppose that $\Lambda_1+\cdots+\Lambda_n=2m$ for some 
positive integer $m$, and let $(a,w_1,\dots,w_m)$ be a solution
of the system of Bethe ansatz equations
\begin{equation}\label{eqBAE2}
{A_+(w_i)}
\prod_{j:j\neq i}
{\theta(w_i-w_j-2\eta)}
=
e^{4a\eta}
{A_-(w_i)}
\prod_{j:j\neq i}{\theta(w_i-w_j+2\eta)}
\qquad i=1,\dots,m,
\end{equation}
such that $w_i\neq w_j\mod\Gamma, (i\neq j)$.
Then $u=\prod Q(x_i)$ with
$Q(x)=e^{ax}\prod_{k=1}^m\theta(x-w_k)$ is a common
eigenfunction of $T(z)$.  
\end{proposition}

\medskip
\noindent{\it Proof:} This is a rephrasing of the first
part of  Prop.\ \ref{p-A4}.
$\square$
\medskip

\medskip
\noindent{\bf Definition.} An eigenfunction of the form
of Prop.\ \ref{p-ba} is called Bethe eigenfunction.

\medskip 

Conversely, let us suppose that $\sum\Lambda_i=2m$, $m\in\Z_{>0}$
and show that all eigenfunctions in a suitable class are of this form. Let,
for a character $\chi:\Gamma\to\C^*$, $\mathcal{H}_m(\chi)$ be the space of
meromorphic functions of $n$ complex variables $x_1,\dots,x_n$ such that
\begin{eqnarray*}
u(\cdots,x_i+1,\cdots)&=&\chi(1)u(\cdots,x_i,\cdots),
\\
u(\cdots,x_i+\tau,\cdots)&=&\chi(\tau)
e^{-\pi im(2x_j+\tau)}u(\cdots,x_i,\cdots),
\end{eqnarray*}
The following result can then easily be
verified using the behavior of the coefficients of the difference
operator $T(z)$.

\begin{lemma} For any character $\chi$ and any $z\in\C$,
$T(z)$ preserves $\mathcal{H}_m(\chi)$.
\end{lemma}

Let for a character $\chi\in\Gamma^*$, $\Sigma(\chi)\subset 
\Theta_m(\chi_0)$ be the set of functions $\epsilon$ so that there exists
a {\em holomorphic} common eigenfunction $u\in \mathcal{H}_m(\chi)$ 
of $T(z)$, $z\in\C$, with eigenvalue $\epsilon(z)$. 

\begin{thm}\label{t-bac}
Suppose that $\Lambda_1+\cdots+\Lambda_n=2m$ for some 
positive integer $m$. Then $\epsilon\in \Sigma(\chi)$ if
and only if there is an eigenfunction $u(x)$  of the form
$u(x)=\prod_{i=1}^nQ(x_i)$ with eigenvalue $\epsilon$, 
such that $Q(x)=e^{ax}\prod_{i=1}^m\theta(
x-w_i)$ for some solution 
$(a,w_1,\dots,w_m)$ of the Bethe ansatz equations
with
\[
\chi(1)=(-1)^me^{a},\qquad 
\chi(\tau)=(-1)^me^{a\tau+2\pi i\sum w_k}.
\]
\end{thm}

\medskip
\noindent{\it Proof:} It remains to show that if $\epsilon\in\Sigma(\chi)$
then this eigenvalue corresponds to a Bethe 
eigenfunction.

Suppose that $\epsilon\in\Sigma(\chi)$
and $u$ is a holomorphic  eigenfunction in $\mathcal{H}_m(\chi)$
with this eigenvalue.
Then, for each $i$, the function $u$ viewed as a 
function of $x_i$ belongs to
$\Theta_m(\chi)$. By setting $z=-x_i$ in the eigenvalue
equation $T(z)u(x)=u(x)$, we see that $u$ is a solution of
the separated equation 
\begin{equation}\label{e-oo}
\prod_{k=1}^n
\theta(x_i+z_k+\eta)\,
T_{x_i}^{-2\eta}u(x)
+
\theta(x_i+z_k-\eta)\,
T_{x_i}^{2\eta}u(x)
=\epsilon(-x_i)u(x).
\end{equation}
Thus this equation admits a non-trivial solution $Q(x_i)$ 
in $\Theta_m(\chi)$ (we consider the remaining variables $x_j$,
$j\neq i$ as fixed).
By the factorization theorem for theta function (see the Appendix),
such a solution is, up to normalization,
 of the form $Q(x_i)=e^{ax}\prod_{j=1}^m\theta(
x_i-w_j)$, and $\chi$ is related to $a$ and $w_k$ by the equations
stated in the Theorem. Setting $x_i=w_j$ in \Ref{e-oo} gives then
the Bethe ansatz equations.
$\square$

\newcommand{\Tirf}{T_{\mathrm{IRF}}}

\section{IRF models with antiperiodic boundary conditions}
We consider here the case where $\Lambda_1=\cdots=\Lambda_n=1$
and show that the commuting transfer matrices $T(z)$ restricted
to functions on $S_0\cap S_1$ are
transfer matrices of IRF models with antiperiodic boundary conditions.
The IRF (interaction-round-a-face)
 models (see \cite{B}, \cite{ABF}, \cite{DJMO})
 of statistical mechanics are two-dimensional lattice models in which
the  configurations over which one sums in the partition function associate
 an element, the height, of a certain set to each pair of neighboring points. The
weight of a configuration is the product of local Boltzmann weights associated
to each ``face'', or square formed by four neighboring points. The
local Boltzmann weight $W(a,b,c,d|z)$ depends on the heights $a,b,c,d$
at the neighboring points and the spectral parameter $z\in\C$. The spectral
parameter can be any fixed number, but more generally, in the
inhomogeneous model, one associates
a spectral parameter to each row and column of faces in the lattice and takes $z$
to be the difference between the row parameter and the column parameter.
In the simplest $sl_2$ case the heights are real numbers
and the local Boltzmann weights are related to
matrix elements of our dynamical $R$-matrix  by
\begin{eqnarray}
 R (z, \lambda = - 2 \eta d) \, e [ c-d ] \otimes e [ b-c] =
   \sum_{a} W(c,b,a,d | z ) \, e [b-a] \otimes e [ a-d] ,
\end{eqnarray} 
where $ c-d,b-c,b-a,a-d \in \{ -1, 1 \}.$ If the latter conditions
are not fulfilled, then $W=0$.\footnote{The normalization chosen here
is so that $W(l,l+1,l+2,l+1|z)=1$ rather than the more common
$\theta(z-2\eta)/\theta(2\eta)$. Also our spectral parameter is 
normalized in a different way than in the literature.}

As shown in \cite{FV4}, the transfer matrix $a(z)+d(z)$ on the
zero weight subspace of the
tensor product of two dimensional representations is then identified
with the row-to-row transfer matrix of the IRF model with periodic
boundary conditions. Here we show that a similar computation applied
to the twisted transfer matrix $b(z)+c(z)$ gives the row-to-row
transfer matrix of the IRF model with antiperiodic boundary conditions.

The formulae are as follows: Let $V(z_1)\otimes\cdot\otimes V(z_n)$ be
a tensor product of two dimensional evaluation modules and assume
that $n$ is odd. Then the $L$-operator
for this module is 
\[
L(z,\lambda)=
R^{(01)} (z - z_1, \lambda - 2 \eta \sum_{i=2}^{n}h^{(i)})
  \cdots
  R^{(0 n-1)} (z-z_{n-1}, \lambda - 2 \eta h^{(n)}) R^{(0n)}(z-z_n,
  \lambda).
\]
 Then the twisted
transfer
matrices, see \ref{s-twisted},
 acting on functions $u(\lambda)$ restricted to the submanifold
given by the equation $\lambda=\eta h$ form a commuting family
parametrized by the spectral parameter. Let
us call these transfer matrices $\Tirf(z)$. The
restriction means that $\Tirf(z)$ preserves functions
$u(\lambda)$ of the form $u(\lambda)=\sum_\mu \delta_{\eta\mu}(\lambda)
u[\mu]$ for some vectors $u[\mu]$ of weight $\mu$. Here $\delta_a(b)$ is
one if $a=b$ and zero otherwise. The restriction to this subspace
of function of $\lambda$ is well-defined if $n$ is odd, as the singularities
of the $L$ operator are at $\lambda\in2\eta\Z$ whereas $\lambda$ takes
values $\eta\mu$, where $\mu$, the weight of a vector in $\otimes V(z_i)$
is an odd number if $n$ is odd. So if $\eta$ is generic, there 
are no singularities at this points. Explicitly,
\begin{eqnarray}
  & \Tirf (z) = \bar b(z)+\bar c(z)=& \nonumber \\ 
  & \sum_{\mu}\mathrm{tr}^{(0)}_{V[\mu]}
 K^{(0)} R^{(01)} (z - z_1, \lambda - 2 \eta \sum_{i=2}^{n}h^{(i)})
   R^{(02)} (z-z_2, \lambda - 2 \eta \sum_{i=3}^{n} h^{(i)}) \cdots  & 
   \nonumber \\
  & R^{(0 n-1)} (z-z_{n-1}, \lambda - 2 \eta h^{(n)}) R^{(0n)}(z-z_n,
  \lambda)T_\lambda^{-2\eta\mu}, & 
\end{eqnarray} 
with
$$  K = \left( \begin{array}{cc} 0 & 1 \\ 1 & 0 
   \end{array} \right), $$ 
The product of $R$-matrices acts on 
 $ V \otimes (\otimes_{i=1}^{n} V) ,$ where the index $(0)$ 
refers to action on the first factor of the tensor product, 
$\mathrm{tr}^{(0)}$
denotes the trace over the first factor. The subsequent 
factors are numbered accordingly from 1 to $n$. \vspace*{0.2cm} \\
The twisted transfer matrices are defined on a $2^n$-dimesnional space.
 This
space has a basis $\delta_{\eta\sum\sigma_i}(\lambda)
e[\sigma_1]\otimes\cdots\otimes e[\sigma_n]$ with $\sigma_i\in\{1,-1\}$.
It is convenient to write this basis in terms of antiperiodic paths:
Let for $a_1,\dots,a_n,a_{n+1}\in\Z+n/2$, such that $|a_i-a_{i+1}|=1$
($i=1,\dots,n$), and $a_{n+1}=-a_1$
\[
|a_1,\dots,a_{n+1}\rangle=\delta_{-2\eta a_{n+1}}(\lambda)
e[a_1-a_2]\otimes\cdots\otimes e[a_n-a_{n+1}]
\]

\begin{proposition} For any antiperiodic path
$|a_1,\dots,a_n,a_{n+1}=-a_1\rangle$, we have
\begin{eqnarray} 
     &  & \Tirf (z) \quad | a_1, \ldots , a_{n+1} \rangle  \qquad  =  \nonumber \\
     &  & \sum_{b_1, \ldots, b_n,b_{n+1}=-b_1} \prod_{i=1}^n W(a_{i+1}, a_i, b_i, b_{i+1})
     \quad | b_1, \ldots , b_{n+1} \rangle \, ,  
\end{eqnarray}
\end{proposition}
\medskip
This expression is, by definition, the transfer matrix of IRF models
with antiperiodic boundary conditions. The partition function for
the IRF model with antiperiodic boundary conditions in one
direction is then $\tr T(w_1)\cdots T(w_m)$. the $w_i$ are the
spectral parameters associated to the $m$ rows, and the $z_i$ are
the spectral parameters associated to the $n$ columns.
\medskip

\noindent{\it Proof of the proposition: } 
Let $e^*[1]$, $e^*[-1]$ be the dual basis to $e[1],e[-1]$.
\begin{eqnarray} 
    &  & \Tirf (z) \quad | a_1, \ldots , a_{n+1} \rangle  \qquad  
=  \nonumber \\
    &  & \sum_\mu \mathrm{tr}^{(0)}_{V[\mu]}
 K^{(0)} R^{(01)} (z-z_1, \lambda - 2 \eta \sum_{i=2}^{n} 
     h^{(i)}) R^{(02)}(z-z_2, \lambda - 2 \eta \sum_{i=3}^{n} h^{(i)}) 
  \nonumber \\
    &  &    \cdots R^{(0n)} (z-z_n, \lambda)
\, \delta_{- 2\eta (a_{n+1}-\mu)} ( \lambda)
     e [a_1 - a_2] \otimes \cdots \otimes e[ a_n - a_{n+1}] =  \nonumber 
\end{eqnarray}
\begin{eqnarray}
    &  & \sum_{\mu} e^{(0) \, \ast} [\mu] K^{(0)} R^{(01)}(z-z_1, \lambda - 2 
     \eta \sum_{i=2}^{n} h^{(i)}) \cdots R^{(0n)} (z-z_n, \lambda) 
     \delta_{-2\eta( a_{n+1}-\mu)} ( \lambda )  \nonumber \\
    &  & \qquad e^{(0)} [\mu] \otimes e [a_1 - a_2] \otimes \cdots \otimes 
      e[ a_n - a_{n+1}]  =  \nonumber \\
    &  & \sum_{b_{n+1}}  e^{(0) \, \ast} [a_{n+1} - b_{n+1}]  K^{(0)} 
      R^{(01)}(z-z_1, \lambda - 2 \eta \sum_{i=2}^{n} h^{(i)}) \cdots R^{(0n)}
      (z-z_n, \lambda)  \nonumber \\
    &  & \qquad \delta_{-2\eta b_{n+1}} ( \lambda)  e^{(0)} 
      [a_{n+1} - b_{n+1}]
      \otimes e [ a_1 - a_2 ] \otimes \cdots \otimes e[a_n - a_{n+1}] = 
      \nonumber  \\
    &  & \qquad \qquad
       \sum_{b_{n+1},b_n}  e^{(0) \, \ast} [a_{n+1} - b_{n+1}]   K^{(0)} 
       R^{(01)}(z-z_1, \lambda - 2 \eta \sum_{i=2}^{n} h^{(i)}) \cdots
       \nonumber \\
    &  & R^{(0 n-1)} (z -z_{n-1}, \lambda - 2 \eta h^{(n)}) \,
        W(a_{n+1},a_n,b_n,b_{n+1}) \delta_{-2\eta b_{n+1}} ( \lambda)
  \nonumber \\
    &  & e^{(0)} [a_n - b_n] \otimes e[a_1
       - a_2] \otimes \cdots \otimes e[a_{n-1} - a_n] \otimes e[b_n - b_{n+1}]
       =   \cdots  =  \nonumber \\
    &  & \sum_{b_{n+1}, \ldots b_i} e^{(0) \, \ast} [a_{n+1}-b_{n+1}]   
      K^{(0)} R^{(01)}(z-z_1, \lambda - 2 \eta \sum_{i=2}^{n} h^{(i)}) 
\cdots
      \nonumber \\
    &  & R^{(0 i-1)} (z-z_{i-1}, \lambda - 2\eta \sum_{j=i}^{n} h^{(j)}) \, 
      \Pi_{j=0}^{n-i} W(a_{n+1-j},a_{n-j},b_{n-j},b_{n+1-j}) 
      \nonumber \\
    &  & \qquad \qquad
         \delta_{-2\eta b_{n+1}} ( \lambda ) e^{(0)} [a_i - b_i] 
         \otimes \nonumber \\ 
    &  & e[ a_1 - a_2]  \otimes 
         e[a_{i-1} - a_{i}] \otimes e[b_{i} - b_{i+1}] \otimes \cdots \otimes 
         e[b_n - b_{n+1}] = \cdots =  \nonumber  \\
    &  & \sum_{b_{n+1}, \ldots, b_1}  e^{(0) \, \ast}  [a_{n+1}-b_{n+1}] 
       K^{(0)} \Pi_{i=1}^{n} W(a_{i+1},a_i,b_i,b_{i+1}) \nonumber \\
    &  & \qquad \delta_{-2\eta b_{n+1}} ( \lambda) e^{(0)} [a_1 
      - b_1] \otimes e[b_1 - b_2] \otimes \cdots \otimes e[b_n -  b_{n+1}] = 
      \nonumber \\
    &  & \sum_{b_{n+1}, \ldots, b_1}  e^{(0) \, \ast}  [a_{n+1}-b_{n+1}] 
        e^{(0)} [- a_1 + b_1] \otimes
      e[b_1 - b_2] \otimes \cdots \otimes e[b_n -  b_{n+1}]   
       \nonumber \\
    &  & \qquad \qquad \delta_{-2\eta b_{n+1}}( \lambda)
        \, \Pi_{i=1}^{n} W(a_{i+1},a_i,b_i,b_{i+1}) \, . \nonumber 
\end{eqnarray} 
By evaluating the linear form $e^\ast[a_{n+1}-b_{n+1}]$ we see that 
only the terms with $b_{n+1}=-b_1$ contribute to the sum. The proof
is complete.
\subsection{Separation of variables for IRF models}
Let us consider the eigenvalue problem for the transfer matrix 
restricted to $S_0\cap S_1$ for $\Lambda_1=\cdots=\Lambda_n=1$ and
$n$ odd. The eigenfunction may be viewed as a function
$u(x_1,\dots,x_n)$ defined for $x_i\in\{-z_i-\eta,-z_i+\eta\}$,
$1\leq i\leq n$. The eigenvalue problem reads
$T(z)u(x)=\epsilon(z)u(x)$, with
\begin{eqnarray}\label{e-sp1}
T(-z)u(x_1,\dots,x_n)&=&\sum_{i=1}^n
\frac{\theta(\lambda-z+x_i)}
{\theta(\lambda)}
\prod_{j\neq i}
\frac{\theta(z-x_j)}
     {\theta(x_i-x_j)}\nonumber\\
&&\bigl(\prod_{k=1}^n\theta(x_i+z_k+\eta\Lambda_k)
u(x_1,\cdots,x_i-2\eta,\cdots,x_n)
\\
&&+\prod_{k=1}^n\theta(x_i+z_k-\eta\Lambda_k)
u(x_1,\dots,x_i+2\eta,\dots,x_n)\bigr)\,,\nonumber
\end{eqnarray}
where $\lambda=-\sum_{k=1}^n(x_k+z_k)$.

\medskip

\begin{thm}\label{t-Spectrum} Suppose that $\Lambda_1=\cdots=\Lambda_n=1$
with $n$ odd, $\eta\not\in\Gamma$,
and that $z_i\neq z_j+2\eta \ell\mod \Gamma$ for $i\neq j$
and $\ell=0,\pm 1$.
Let $T(z)=b(z)+c(z)$ be the transfer matrix restricted
to the $2^n$ dimensional space of functions on $S_0\cap S_1$. Then 
a function $\epsilon(z)$ is a common eigenvalue of the
transfer matrices $T(z)$, $z\in\C$, if and only if
\begin{enumerate}
\item[(i)]
$\epsilon\in\Theta_n(\chi)$
with $\chi(1)=(-1)^n$, $\chi(\tau)=(-1)^ne^{2\pi i\sum z_j}$
and 
\item[(ii)]
$\epsilon$ obeys the quadratic relations
\begin{equation}\label{e-qr}
  \epsilon(z_i-\eta)\epsilon(z_i+\eta)=\prod_{k=1}^n
\theta(z_k-z_i+2\eta)
\theta(z_k-z_i-2\eta), \qquad i=1,\dots,n.
\end{equation}
\end{enumerate}
\end{thm}

To prove this theorem, let $S=\times_{i=1}^n\{-z_i-\eta,-z_i+\eta\}$.
Let us first assume that $u(x)$, $x\in S$
is a common eigenfunction of $T(z)$ with eigenvalue $\epsilon(z)$. 
In particular $u$ does not vanish identically on $S$. From the
transformation properties of $T(z)$ under shifts of $z$ by $\Gamma$
we see that $\epsilon(z)$ has to belong to $\Theta_n(\chi)$.
Then setting $z=-x_i$ in the eigenvector equation $T(z)u(x)=
\epsilon(z)u(x)$, we get the separated equation
\begin{equation}\label{e-ooo}
\prod_{k=1}^n
\theta(x_i+z_k+\eta)\,
T_{x_i}^{-2\eta}u(x)
+
\theta(x_i+z_k-\eta)\,
T_{x_i}^{2\eta}u(x)
=\epsilon(-x_i)u(x).
\end{equation}
Inserting the two values of $x_i$ yields 
\begin{equation}\label{e-ttt}
\begin{array}{c}
\prod_{k=1}^n
\theta(z_k-z_i+2\eta)
\,
u(x_1,\dots,-z_i-\eta,\dots,x_n)
=
\epsilon(z_i-\eta)
u(x_1,\dots,-z_i+\eta,\dots,x_n),
\\
 \\
\prod_{k=1}^n
\theta(z_k-z_i-2\eta)
\,
u(x_1,\dots,-z_i+\eta,\dots,x_n)
=
\epsilon(z_i+\eta)
u(x_1,\dots,-z_i-\eta,\dots,x_n).
\end{array}
\end{equation}
By multiplying these two equations we see that $\epsilon$ must
obey the identity \Ref{e-qr}, provided we can prove that,
for at least one choice of $x_j$, the product 
\[
u(x_1,\dots,-z_i-\eta,\dots,x_n)\,u(x_1,\dots,-z_i+\eta,\dots,x_n)
\]
does not vanish. This follows from the fact that $u$ is not identically
zero, so that at least one factor of this product is nonzero, and
that 
the product of theta functions on the left-hand side of \Ref{e-ttt}
is not zero with our assumption on the $z_j$, so that also the other
factor is nonzero.

We have thus  shown that a necessary condition for a funtion $\epsilon(z)$
to be a common eigenvalue is that $\epsilon$ is a theta funtion obeying
the quadratic relations \Ref{e-qr}.

Conversely, let us suppose that $\epsilon$ obeys \Ref{e-qr}. Then,
for every $i$, the system of equations
\[
\prod_{k=1}^n
\theta(x+z_k+\eta)\,
Q(x-2\eta)
+
\theta(x+z_k-\eta)\,
Q(x+2\eta)
=\epsilon(-x)Q(x), 
\quad x=-z_i\pm\eta,
\]
admits a non-trivial solution $Q_i(x)$. It follows that
$u(x)=\prod_{i=1}^nQ_i(x_i)$ obeys \Ref{e-ooo}. Thus, for any
$x\in S$,
$(T(z)-\epsilon(z))u(x)$, viewed as a function of $z$ is a
theta function in $\Theta_n(\chi)$ vanishing at $n$ distinct points
$-x_1,\dots,-x_n$. As explained in the Appendix, this implies that
that either the function vanishes identically or 
$-\sum x_i=\sum z_i$. Since for $x\in S$, $\sum x_i$ is $-\sum z_i$ plus
an {\em odd} multiple of $\eta$, the latter alternative cannot hold.
Thus $u$ is an eigenfuntion with eigenvalue $\epsilon$.
The Theorem is proven.

\medskip

\section{Conclusions}

In this paper we have studied integrable models
associated to elliptic curves by Sklyanin's method of 
separation of variables. We have shown how in the elliptic version
of the Gaudin model, previously considered in \cite{EFR}, this method
implies the completeness of Bethe eigenfunction in Verma modules,
in the sense that for every eigenvalue one has a Bethe eigenfunction
(although in the case of degenerate eigenvalues one cannot exclude
the existence of additional eigenfunction not of Bethe type).

In the difference case, the analogue of the transfer matrix in 
the separated variables was shown to be a twisted transfer matrix
associated to a representation of $E_{\tau,\eta}(sl_2)$ by
difference operators. The eigenvalue problem can be posed in two
different cases. In the continuous case, the variables are assumed
to be complex and one asks for eigenfunctions with theta functions
properties, as one does in the differential case. These eigenfunctions
are of Bethe ansatz type. In the discrete case, the variables are
assumed to take values in a finite set. The eigenvalue problem can
still be solved by separating variables yielding eigenfunctions
of the (inhomogeneous)
row-to-row transfer matrix of IRF models with antiperiodic boundary
conditions. In both cases one has a completeness result.

Let us conclude by mentioning some open questions. In the
difference case one has the transfer matrix in the ``separated
variables''. Is there a ``quantum Radon transform'' as in 
the differential case, which maps the transfer matrix of the
IRF model with periodic boundary conditions to this transfer 
matrix? Is there an analogue of the equivalence between the
local problem and the global problem of \cite{TV} in the
elliptic case?

\appendix
\section{Elliptic polynomials}\label{s-A}
\subsection{Theta functions}
Let $\mathrm{Im}\,\tau>0$ and set $q=e^{2\pi i\tau}$.
Let $\Gamma$ be the lattice $\Z+\tau\Z$ and $\Gamma^*\simeq (\C^\times)^2$ 
the
group of  group homomorphisms $\Gamma\to\C^\times$. 
Let $\phi$ be the homomorphism
$\phi:\chi\mapsto \frac1{2\pi i}(\ln\,\chi(\tau)
-\tau\ln\,\chi(1))$ (or, more invariantly, $\frac1{2\pi i}
(\omega_1\ln\,\chi(\omega_2)
-\omega_2\ln\,\chi(\omega_1))$, for any oriented basis $\omega_1,\omega_2$
of $\Gamma$)
from $\Gamma^*$ to the elliptic curve $E=\C/\Gamma$.

For $\chi\in\Gamma^*$ let
$\Theta_k(\chi)$ be the space of theta functions of level
$k$ and character $\chi$. It consists of
entire holomorphic functions $f(z)$ such that
$f(z+r+s\tau)=\chi(r+s\tau)\exp(-\pi i  k(s^2\tau+2sz))f(z)$
for all $r+s\tau\in\Gamma$.

The dimension of $\Theta_k(\chi)$ is zero if $k<0$. It is $k$ if $k\geq1$.
The dimension of $\Theta_0(\chi)$ is one if $\phi(\chi)=0$ and vanishes 
otherwise.

We have the unique factorization result:
\begin{proposition}\label{unique factorization}
The function of $z\in\C$
\[
f(a,w;z)=e^{az}\prod_{j=1}^k\theta(z-w_j)
\]
belongs to $\Theta_k(\chi)$, with
$\chi(r+s\tau)=(-1)^{(r+s)k}
e^{ra+s(a\tau+2\pi i\sum_j w_j)}$. Every function in $\Theta_k(\chi)$
is of the form $C\cdot f(a,w;z)$ for some constant $C$ and this
representation is unique up to permutation of the $w_j$
if one requires 
 the $w_j$ to be in the fundamental domain
$F=\{x+y\tau\,|\,x,y\in[0,1)\}$.\end{proposition}
\begin{prf}
It follows from the transformation properties of theta functions
that the number of zeros $(2\pi i)^{-1}\int_{\partial F}d\ln g$,
counted with multiplicities,  in
$F$ of
a theta function $g\in\Theta_k(\chi)$ is $k$. If $w_1,\dots,w_k$
denote the zeros of $g$ then $g(z)/f(a,w;z)$ is doubly periodic
and regular, thus constant. Uniqueness follows from the fact that
$a$ is uniquely determined by the zeros $w_j$ and the character 
$\chi$.
\end{prf}
\begin{corollary}\label{cor1}
Let $E$ be the elliptic curve $\C/\Gamma$, and, for $k=1, 2,\dots$,
$S^k(E)=E/S_k$ be
its $k\th$ symmetric power.
The map $\mathbb{P}(\Theta_k(\chi))\to S^k(E)$ sending a function to
the set of its zeros modulo $\Gamma$, is injective. Its image
consists of classes $[w_1,\dots,w_k]$ such that $\sum_jw_j
=\phi(\chi)+k\delta$, where $\delta$ is the image in $E$ of $(1+\tau)/2$.
\end{corollary}
\subsection{Interpolation}
\begin{thm}\label{t-A3}
Suppose $z_1,\dots,z_k\in\C$ are pairwise distinct modulo $\Gamma$ and
$\chi\in\Gamma^*$ is such that $\sum_{i=1}^k z_i\neq \phi(\chi)+k\delta
\mod \Gamma$. 
Then for any $y_1,\dots,y_k\in\C$ 
there exists a unique function $f\in\Theta_k(\chi)$ such that
$f(z_i)=y_i$ for all $i=1,\dots,k$.
\end{thm}

The interpolation formula giving $f$ is
\[
f(z)=\sum_{j=1}^my_je^{2\pi ia(z-z_j)}
\frac{\theta(z-z_j+b)}{\theta(b)} 
\prod_{l:l\neq j}\frac{\theta(z-z_l)}{\theta(z_j-z_l)}.
\]
Here $a$ and $b$ are such that the character of all terms in the
sum are equal to $\chi$, so
\begin{eqnarray*}
a&=&\frac1{2\pi i}\ln\chi(1)-\frac k2,\\
b&=&\frac1{2\pi i}\left(\tau\ln\chi(1)-\ln\chi(\tau)\right)
+\sum_{j=1}^kz_j-k\frac{1+\tau}2\,,
\end{eqnarray*}
for any choice of the branch of the logarithm.
The assumption on $\chi$ ensures that $b\not\in\Gamma$ so that
the denominator $\theta(b)$ does not vanish.

The function is unique since the difference of any two is a theta
function vanishing
at $m$ points. By Corollary \ref{cor1}, with our assumption on $\chi$,
it must vanish identically.

\subsection{Difference equations}
Consider the following problem arising in integrable
models. Given $\gamma\in\C-\Gamma$ and functions $A_+(z)\in
\Theta_k(\chi_+)$,
$A_-(z)\in\Theta_k(\chi_-)$, find $\epsilon(z)$ such that
the difference equation 
\begin{equation}\label{eqDE}
A_+(z)Q(z-\gamma)
+
A_-(z)Q(z+\gamma)
=
\epsilon(z)Q(z),
\end{equation}
has a non-trivial solution $Q(z)$ in some $\Theta_m(\chi)$.

A necessary condition is that all terms are theta functions with
the same character. So $\epsilon$ has to be of level $k$ and
$\chi_+(1)=\chi_-(1)$, 
$\chi_+(\tau)e^{2\pi i\gamma m}=
\chi_-(\tau)e^{-2\pi i\gamma m}=$ 
character of $\epsilon$.

A pair $(\epsilon,Q)$ of theta functions (with $Q$ non-trivial)
obeying the difference
equation \Ref{eqDE} will be called an elliptic polynomial solution
of \Ref{eqDE}.

\begin{proposition}\label{p-A4}
Suppose $A_\pm(z)\in\Theta_k(\chi_\pm)$ with 
$\chi_+(r+s\tau)=\chi_-(r+s\tau)e^{-4\pi i\gamma ms}$ for all
$r+s\tau\in\Gamma$ and some positive integer $m$.
If $a,w_1,\dots,w_m$ are solutions of the system of equations
\begin{equation}\label{eqBAE}
{A_+(w_i)}
\prod_{j:j\neq i}
{\theta(w_i-w_j-\gamma)}
=
e^{2a\gamma}
{A_-(w_i)}
\prod_{j:j\neq i}{\theta(w_i-w_j+\gamma)}
\qquad i=1,\dots,m.
\end{equation}
such that $w_i\neq w_j \mod\Gamma$ for all $i\neq j$,
then the functions
\begin{equation}\label{eqQ}
Q(z)=e^{az}\prod_{j=1}^m\theta(z-w_j)
\end{equation}
and
\begin{equation}\label{eqE}
\epsilon(z)=\frac{A_+(z)Q(z-\gamma)+A_-(z)Q(z+\gamma)}
{Q(z)}
\end{equation}
form an elliptic polynomial solution of \Ref{eqDE}. 
Conversely, if $(\epsilon,Q)$ is an elliptic polynomial solution
of \Ref{eqDE}, then there exists a solution $a,w_1,\dots,w_m$
of the system \Ref{eqBAE} such that  $Q$ is of the form \Ref{eqQ}
(up to a multiplicative constant) and  $\epsilon$ is given by
\Ref{eqE}.
\end{proposition}

\begin{prf} Let $Q$ be the function defined 
in \Ref{eqQ}. The ratio \Ref{eqE} obeys the transformation
property of a theta function, but may be singular at the
zeros $w_i$ of $Q$.

The $i\th$ equation in the system, or more precisely, the
equation equivalent to it if $\gamma\in\C-\Gamma$:
\[
A_+(w_i)e^{-\gamma a}\prod_{j=1}^m\theta(w_i-w_j-\gamma)
+
A_-(w_i)e^{\gamma a}\prod_{j=1}^m\theta(w_i-w_j+\gamma)
=0,
\]
is the condition that the left-hand side of
the \Ref{eqDE} vanishes at $w_i$. Thus, if $a, w_1,\dots,w_m$
obey the system of equations, then the
quotient $\epsilon$ of this left-hand side by $Q(z)$ is regular
at $w_i$ and thus everywhere
and $(\epsilon,Q)$ is an elliptic polynomial solution.

Suppose now that $(\epsilon,Q)$ is an  elliptic polynomial
solution. By Proposition \ref{unique factorization},
$Q$ can be written in the form \Ref{eqQ} if we normalize it
properly. The points $w_i$ are
the zeros (modulo $\Gamma$) of $Q$. Then the 
left-hand side of the difference equation vanishes at $w_i$ so
the $w_i$ are a solution of \Ref{eqBAE}
\end{prf}

\medskip
\noindent{\bf Remark.} The system of equations is usually written
in the form
\[
\prod_{j:j\neq i}
\frac
{\theta(w_i-w_j-\gamma)}
{\theta(w_i-w_j+\gamma)}
=
e^{2a\gamma}
\frac
{A_-(w_i)}
{A_+(w_i)},
\qquad i=1,\dots,m.
\]


\begin{thebibliography}{ccc}
\bibitem{ABF}G. E. Andrews, R. J. Baxter, P. J. Forrester,
{\em Eight-vertex {\rm SOS} model and generalized Rogers--%
Ramanujan-type identities}. J. Stat. Phys. 35 (1984), no. 3-4, 193--266
\bibitem{B}   R. J. Baxter, {\it Exactly solved models of statistical
mechanics}, Academic Press, New York 1982
\bibitem{DJMO} E. Date, M. Jimbo, T. Miwa, M. Okado,
{\em Fusion of the eight vertex SOS model},
Lett. Math. Phys. 12 (1986), no. 3, 209--215
\bibitem{EFR} B. Enriquez, B. Feigin and V. Rubtsov,
{\it Separation of variables for Gaudin-Calogero systems},
{\em Comp.\ Math.} { 110}
(1998),
1--16.
\bibitem{FV1} G. Felder and A. Varchenko, 
{\it Integral representation of solutions of the 
elliptic Knizhnik-Zamolodchikov-Bernard equation},
Int. Math. Res. Notices, No.\ 5(1995), 221--233
\bibitem{FV2} G. Felder and A. Varchenko, 
{\it Three formulae for eigenfunctions of integrable
Schr\"odinger operators}, Comp.\ Math.\ 107 (1997), 143--175
\bibitem{FV3} G. Felder and A. Varchenko, 
{\it On representation theory of the elliptic quantum
group $E_{\tau,\eta}(sl_2)$}, 
Commun.\ Math.\ Phys.\  181 (1996), 741--761
\bibitem{FV4} G. Felder and A. Varchenko, 
{\it Algebraic Bethe anstaz for the elliptic quantum 
group $E_{\tau,\eta}(sl_2)$}, 
Nucl.\ Phys. B 480 (1996), 485--503
\bibitem{F}   E. Frenkel,
{\it Affine algebras, Langlands duality and Bethe ansatz}, 
Proceedings of the XIth International Congress of Mathematical Physics
(Paris, 1994),
606--642,
Internat. Press, Cambridge, MA,
1995
\bibitem{S}   E. Sklyanin,{\em Separation of variables for the Gaudin
model}, J. Sov. Math. 
{ 47}
(1989),
2473--2488;
{\it
Separation of variables. New trends},
in {\em Quantum field theory, integrable models and beyond}
(Kyoto, 1994). 
Progr. Theoret. Phys. Suppl.
{ 118}
(1995), 
35--60.
\bibitem{ST}  E. Sklyanin and T. Takebe,
{\it Separation of varibales in the elliptic Gaudin model},
solv-int/9807008
\bibitem{TV}  V. Tarasov and A. Varchenko,
{\it Completeness of Bethe vectors and difference equations with
 regular singular points},
Internat. Math. Res. Notices 1995, no. 13, 637--669
\bibitem{WW} 
E. Whittaker and G. Watson,
{\it A Course of Modern Analysis}, 4th ed.,
Cambridge University Press,
1927
\end{thebibliography}
\end{document}